\documentclass[journal]{IEEEtran}
\ifCLASSINFOpdf
\usepackage{cite}
\usepackage{hyperref}
\usepackage[pdftex]{graphicx}
\usepackage[export]{adjustbox}
\usepackage{kantlipsum}

\graphicspath{{figures/}}
\else
\fi
\usepackage{dblfloatfix} 
\usepackage{kantlipsum} 

\usepackage{amsmath}
\usepackage{amssymb}
\usepackage{amsthm}
\usepackage{mathtools}
\usepackage{stackengine}

\DeclareMathOperator*{\argmin}{arg\,min}

\usepackage{algorithmic, algorithm}

\usepackage{booktabs}

\usepackage{calc}
\usepackage{xcolor}
\usepackage{dsfont}
\label{comments}

\numberwithin{theorem}{section} 

\usepackage{booktabs}
\usepackage{multirow}
\definecolor{myBlue}{RGB}{219, 48, 122}

\usepackage{pgf,tikz,pgfplots}
\usetikzlibrary{arrows,shapes}
\usetikzlibrary{plotmarks}
\usepgfplotslibrary{groupplots}
\pgfplotsset{compat=newest}
\pgfplotsset{plot coordinates/math parser=false}
\usetikzlibrary{positioning,spy}
\usetikzlibrary{backgrounds}
\usepackage{comment}

\definecolor{my_yellow}{rgb}{0.9,0.8,0.05}
\definecolor{my_green}{rgb}{0.2,0.5,0.1}

\begin{document}

	\title{Deep Coded Aperture Design: An End-to-End Approach for Computational Imaging Tasks}
	\author{Jorge Bacca,~\IEEEmembership{Student~Member,~IEEE,} and Tatiana Gelvez,~\IEEEmembership{Student~Member,~IEEE} and Henry Arguello,~\IEEEmembership{Senior~Member,~IEEE}
	\thanks{J. Bacca and H. Arguello are with the Department of Computer Science,
Universidad Industrial de Santander, Bucaramanga, Colombia, 680002 e-mail:
jorge.bacca1@correo.uis.edu.co; henarfu@uis.edu.co}
	\thanks{T. Gelvez is with the Department
of Electrical Engineering, Universidad Industrial de Santander, Bucaramanga, Colombia, 680002, e-mail: tatiana.gelvez@correo.uis.edu.co.}\vspace{-1em}
\thanks{This work was supported by the VIE of Universidad Industrial de Santander through “Cámara difractiva compacta y de bajo costo para la adquisición comprimida de imágenes espectrales en la industria agrícola Santandereana.” under  Project 2699. The work of Tatiana was supported by the Academy of Finland with project no. 318083 }}


	\maketitle
	
	\begin{abstract}
	
	Covering from photography to depth and spectral estimation, diverse computational imaging (CI) applications benefit from the versatile modulation of coded apertures (CAs). The light wave fields as space, time, or spectral can be modulated to obtain projected encoded information at the sensor that is then decoded by efficient methods, such as the modern deep learning decoders. Despite the CA can be fabricated to produce an analog modulation, a binary CA is mostly preferred since easier calibration, higher speed, and lower storage are achieved. As the performance of the decoder mainly depends on the structure of the CA, several works optimize the CA ensembles by customizing regularizers for a particular application without considering critical physical constraints of the CAs. This work presents an end-to-end (E2E) deep learning-based optimization of CAs for CI tasks. The CA design method aims to cover a wide range of CI problems easily changing the loss function of the deep approach. The designed loss function includes regularizers to fulfill the widely used sensing requirements of the CI applications. Mainly, the regularizers can be selected to optimize the transmittance, the compression ratio and the correlation between measurements, while a binary CA solution is encouraged, and the performance of the CI task is maximized in applications such as restoration, classification, and semantic segmentation. 

	\end{abstract}
	
	

	\IEEEpeerreviewmaketitle

	\section{Introduction}

	
	\IEEEPARstart{T}{he} 	\textcolor{black}{light wavefront emitted or reflected by an object is modeled mathematically with the plenoptic function, considering that some application fields of the light as, spectral~\cite{arguello2014colored}, 
	temporal~\cite{baraniuk2017compressive}, spatial~\cite{llull2013coded}, polarization~\cite{tsai2013coded},
	depth~\cite{zhou2009coded},
	amplitude~\cite{zhang2017hadamard}, 
	phase~\cite{bacca2019super},
	 and angular views~\cite{marwah2013compressive},
	 can be encoded by an optical element known as coded aperture (CA).} This modulation is useful in applications where directly capturing the incoming light is impractical or infeasible~\cite{tsai2015spectral,lai2018coded,radwell2014single}. Besides, the CA alleviates the dimensionality mismatch without scanning through each of the plenoptic \textcolor{black}{fields}, since a proper encoding reduces the uncertainty in the recovery of specific wavefront dimensions through computational processing after acquiring projected encoded measurements with an optical coding system~\cite{cao2016computational}; we refer to as optical coding system to any architecture that incorporates a CA in the setup. 

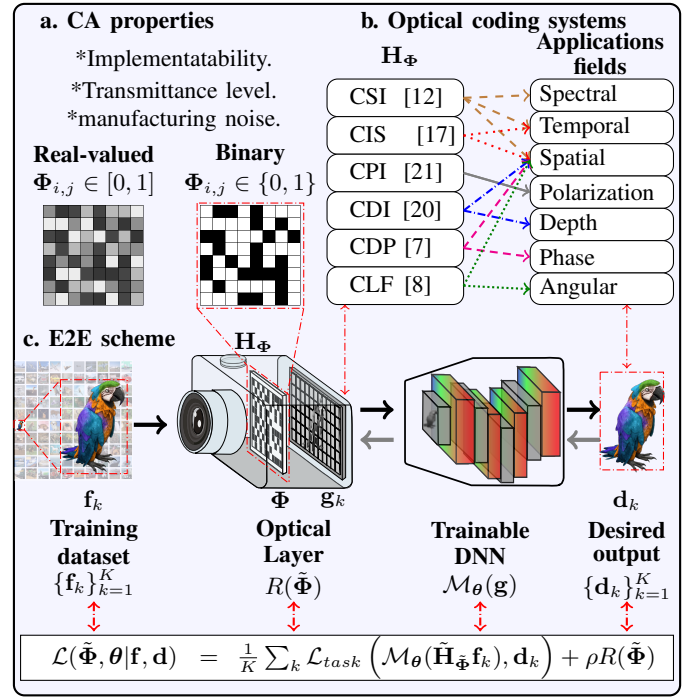
\begin{figure}[!b]
	\begin{tikzpicture}
 \tikzstyle{every node}=[font=\small]
 
\filldraw [fill=blue!4, draw=black , rounded corners, thick] (0,-0) rectangle (1\linewidth,-9.2); 

\node[inner sep=0pt] (binaryCA) at (0.12\linewidth,-3.35)
 {\includegraphics[width=.15\linewidth]{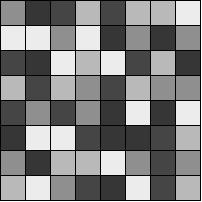}};
\node[inner sep=0pt] (binary2CA) at (0.355\linewidth,-3.35)
 {\includegraphics[width=.15\linewidth]{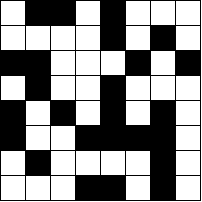}};


\node at (0.24\linewidth,-0.75) {*Implementatability.};
\node at (0.24\linewidth,-1.15) {*Transmittance level.};
\node at (0.24\linewidth,-1.55) {*manufacturing noise.};
 
\node at (0.17\linewidth,-0.245cm) {\textbf{a. CA properties}};
\node at (0.72\linewidth,-0.245cm) {\textbf{b. Optical coding systems}};
\node at (0.12\linewidth,-2.0cm) {\textbf{Real-valued}};
\node at (0.355\linewidth,-2.0cm) {\textbf{Binary}};
\node at (0.12\linewidth,-2.4cm) {$\boldsymbol{\Phi}_{i,j} \in [0,1]$};
\node at (0.355\linewidth,-2.4cm) {{$\boldsymbol{\Phi}_{i,j} \in \{0,1\}$}};

 
\filldraw [fill=white, draw=black, rounded corners] (0.47\linewidth,-1.0) rectangle (0.675\linewidth,-1.5) node[xshift=-0.095\linewidth, yshift=0.8cm] {$\mathbf{H}_{\boldsymbol{\Phi}}$};
\filldraw [fill=white, draw=black,rounded corners] (0.47\linewidth,-1.5) rectangle (0.675\linewidth,-2.0);
\filldraw [fill=white, draw=black, rounded corners] (0.47\linewidth,-2.0) rectangle (0.675\linewidth,-2.5);
\filldraw [fill=white, draw=black, rounded corners] (0.47\linewidth,-2.5) rectangle (0.675\linewidth,-3.0);
\filldraw [fill=white, draw=black, rounded corners] (0.47\linewidth,-3.0) rectangle (0.675\linewidth,-3.5);
\filldraw [fill=white, draw=black,, rounded corners] (0.47\linewidth,-3.5) rectangle (0.675\linewidth,-4.0);

\node at (0.49\linewidth,-3.75) [anchor=west] {CLF \hspace{0.1em}\cite{marwah2013compressive}};
\node at (0.49\linewidth,-3.25) [anchor=west] {CDP~\cite{bacca2019super}};
\node at (0.49\linewidth,-2.75) [anchor=west] {CDI \hspace{0.15em}\cite{levin2007image}};
\node at (0.49\linewidth,-2.25) [anchor=west] {CPI \hspace{0.35em}\cite{tu2017division}};
\node at (0.49\linewidth,-1.75) [anchor=west] {CIS \hspace{0.9em}\cite{duarte2008single}};
\node at (0.49\linewidth,-1.25) [anchor=west] {CSI \hspace{0.35em}\cite{cao2016computational}};

\filldraw [fill=white, draw=black, rounded corners] (0.7753\linewidth,-1.0) rectangle (0.99\linewidth,-1.43) node[xshift=-0.11\linewidth, yshift=0.9cm] {\textbf{\textcolor{black}{Applications}}}
node[xshift=-0.11\linewidth, yshift=0.6cm] {\textbf{\textcolor{black}{fields}}};
\filldraw [fill=white, draw=black, rounded corners] (0.775\linewidth,-1.43) rectangle (0.99\linewidth,-1.86);
\filldraw [fill=white, draw=black, rounded corners] (0.775\linewidth,-1.86) rectangle (0.99\linewidth,-2.29);
\filldraw [fill=white, draw=black,rounded corners] (0.775\linewidth,-2.29) rectangle (0.99\linewidth,-2.72);
\filldraw [fill=white, draw=black,rounded corners] (0.775\linewidth,-2.72) rectangle (0.99\linewidth,-3.15);
\filldraw [fill=white, draw=black,rounded corners] (0.775\linewidth,-3.15) rectangle (0.99\linewidth,-3.58);
\filldraw [fill=white, draw=black,rounded corners]
(0.775\linewidth,-3.58) rectangle (0.99\linewidth,-4.00);

\node at (0.775\linewidth,-3.785) [anchor=west] {Angular};
\node at (0.775\linewidth,-3.365) [anchor=west] {Phase};
\node at (0.775\linewidth,-2.935) [right] {Depth};
\node at (0.775\linewidth,-2.505) [right] {Polarization};
\node at (0.775\linewidth,-2.075) [right] {Spatial};
\node at (0.775\linewidth,-1.645) [anchor=west] {Temporal};
\node at (0.775\linewidth,-1.215) [anchor=west] {Spectral};

\draw [<->, red, densely dashdotted] (4.4, -5.2) -- (4.4, -4.0);

\draw [->,brown, dashed, thick] (0.675\linewidth,-1.25) -- (0.775\linewidth,-1.645);
\draw [->,brown, dashed, thick] (0.675\linewidth,-1.25) -- (0.775\linewidth,-1.215);
\draw [->,brown, dashed, thick] (0.675\linewidth,-1.25) -- (0.775\linewidth,-2.075);
\draw [->,red, dotted, thick] (0.675\linewidth,-1.75) -- (0.775\linewidth,-1.645);
\draw [->,red, dotted, thick] (0.675\linewidth,-1.75) -- (0.775\linewidth,-2.075);
\draw [->, gray,solid, thick] (0.675\linewidth,-2.25) -- (0.775\linewidth,-2.505);
\draw [->, gray, solid, thick] (0.675\linewidth,-2.25) -- (0.775\linewidth,-2.505);
\draw [->,blue,densely dashdotted, thick] (0.675\linewidth,-2.75) -- (0.775\linewidth,-2.075);
\draw [->,blue,densely dashdotted, thick] (0.675\linewidth,-2.75) -- (0.775\linewidth,-2.935);
\draw [->,magenta, densely dashed, thick] (0.675\linewidth,-3.25) -- (0.775\linewidth,-3.365);
\draw [->,magenta, densely dashed, thick] (0.675\linewidth,-3.25) -- (0.775\linewidth,-2.075);
\draw [->,black!50!green!100, densely dotted, thick] (0.675\linewidth,-3.75) -- (0.775\linewidth,-3.795);
\draw [->,black!50!green!100, densely dotted, thick] (0.675\linewidth,-3.75) -- (0.775\linewidth,-2.075);


 	\node[inner sep=0pt] (system) at (0.5\linewidth,-5.55)
 {\includegraphics[width=1\linewidth]{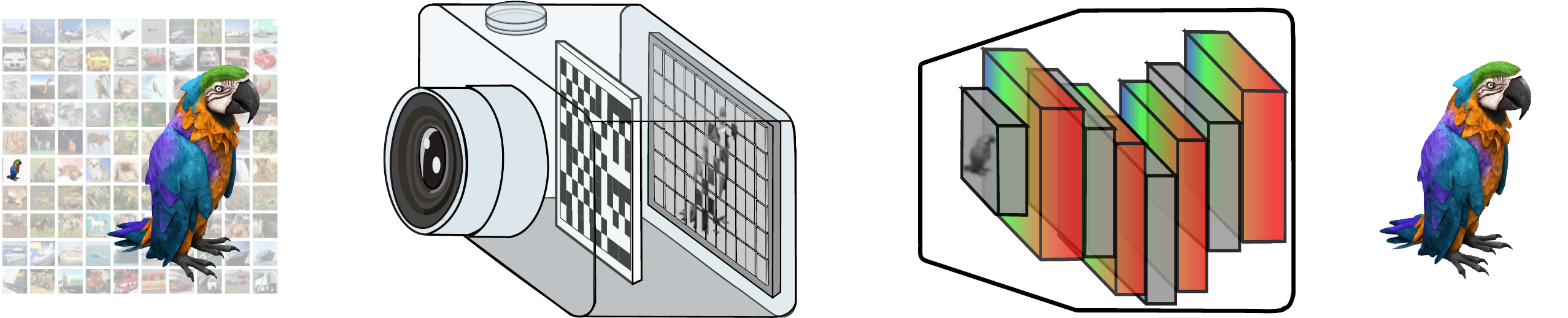}};

\node[anchor=west] at (0.0\linewidth,-4.45cm) {\textbf{c. E2E scheme}};

	\draw[red,densely dashdotted] (0.275\linewidth,-2.65) rectangle (0.435\linewidth,-4.05);
	
	\draw[red,densely dashdotted] (0.275\linewidth,-4.05) -- (0.35\linewidth,-4.81cm);
	\draw[red,densely dashdotted] (0.435\linewidth,-4.05) -- (0.41\linewidth,-5.2cm);

	\draw[red,densely dashdotted] (0.35\linewidth,-4.81) -- (0.35\linewidth,-5.95);
	\draw[red,densely dashdotted] (0.35\linewidth,-5.95) -- (0.413\linewidth,-6.38);
	\draw[red,densely dashdotted] (0.413\linewidth,-5.2) -- (0.413\linewidth,-6.36);
	\draw[red,densely dashdotted] (0.354\linewidth,-4.81) -- (0.413\linewidth,-5.2);
\draw[red,densely dashdotted] (0.0\linewidth,-5.54) rectangle (0.016\linewidth,-5.67);
\draw[red,densely dashdotted] (0.07\linewidth,-5.0) rectangle (0.17\linewidth,-6.25);
\draw[red,densely dashdotted] (0.016\linewidth,-5.54) -- (0.07\linewidth,-5.0cm);
\draw[red,densely dashdotted] (0.016\linewidth,-5.67) -- (0.07\linewidth,-6.25cm);

\draw[<->,red, densely dashdotted, thick] (0.12\linewidth,-7.95) -- (0.12\linewidth,-8.35);
\draw[<->,red, densely dashdotted, thick] (0.42\linewidth,-7.95) -- (0.42\linewidth,-8.35);
\draw[<->,red, densely dashdotted, thick] (0.7\linewidth,-7.95) -- (0.7\linewidth,-8.35);
\draw[<->,red, densely dashdotted, thick] (0.92\linewidth,-7.95) -- (0.92\linewidth,-8.35);

\draw[->,black, solid, ultra thick] (0.18\linewidth,-5.6) -- (0.23\linewidth,-5.6); 
\draw[->,black, solid, ultra thick] (0.52\linewidth,-5.4) -- (0.57\linewidth,-5.4); 
\draw[->,black, solid, ultra thick] (0.83\linewidth,-5.4) -- (0.88\linewidth,-5.4); 
\draw[<-,black!50, solid, ultra thick] (0.83\linewidth,-5.8) -- (0.88\linewidth,-5.8); 
\draw[<-,black!50, solid, ultra thick] (0.52\linewidth,-5.8) -- (0.57\linewidth,-5.8); 

\node at (0.12\linewidth,-7.0) {\textbf{Training}};
\node at (0.12\linewidth,-7.35) {\textbf{dataset}};
\node at (0.12\linewidth,-7.7) {$\{\mathbf{f}_k\}_{k=1}^K$};
\node at (0.42\linewidth,-7.0) {\textbf{Optical}};
\node at (0.42\linewidth,-7.35) {\textbf{Layer}};
\node at (0.42\linewidth,-7.75cm) {$ R(\tilde{\boldsymbol{\Phi}})$};
\node at(0.7\linewidth,-7.0) {\textbf{Trainable}};
\node at(0.7\linewidth,-7.35) {\textbf{DNN}};
\node at (0.7\linewidth,-7.75) {$\mathcal{M}_{\boldsymbol{\theta}}(\mathbf{g})$};
\node at(0.92\linewidth,-7.0) {\textbf{Desired}};
\node at(0.92\linewidth,-7.35) {\textbf{output}};
\node at (0.92\linewidth,-7.75) {$\{\mathbf{d}_k\}_{k=1}^K$};
\node at (0.92\linewidth,-6.6) {$\mathbf{d}_k$};
\node at (0.48\linewidth,-6.6) {$\mathbf{g}_k$};
\node at (0.12\linewidth,-6.6) {$\mathbf{f}_k$};
\node at (0.4\linewidth,-6.6cm) {$\boldsymbol{\Phi}$};
\node at (0.36\linewidth,-4.5cm) {$\mathbf{H}_{\boldsymbol{\Phi}}$};

\draw [<->, red, densely dashdotted] (0.92\linewidth, -4.0) -- (0.92\linewidth, -4.85);

\draw[red,densely dashdotted] (0.88\linewidth,-4.9) rectangle (0.975\linewidth,-6.2);

\filldraw [fill=white, draw=black] (0.01\linewidth,-8.4) rectangle (0.99\linewidth,-9) node[xshift=-0.48\linewidth, yshift=0.3cm] {$\mathcal{L}(\tilde{\boldsymbol{\Phi}}, \boldsymbol{\theta} | \mathbf{f}, \mathbf{d}) \hspace{2mm} = \hspace{1.5mm} \frac{1}{K} \sum_{k} \mathcal{L}_
	{task}\left(\mathcal{M}_{\boldsymbol{\theta}}(\tilde{\mathbf{H}}_{\tilde{\boldsymbol{\Phi}}}\mathbf{f}_k ), \mathbf{d}_k\right)+ \rho R(\tilde{\boldsymbol{\Phi}})$};
\end{tikzpicture}\vspace{-1.5em}
	\caption{Proposed E2E Approach. (\romannumeral 1) The sensing protocol is modeled as a learnable optical layer whose trainable parameter is the CA. A set of scenes passes through the optical layer to obtain the projected measurements that enter to the hidden convolutional layers up to the loss function to achieve the specific task. The error propagates back up to the CA.}
	\label{fig:proposed_method}
	\end{figure}

	Analytically, the CA is modeled as a tensor array where each spatial location has a particular response to the incoming wavefront. \textcolor{black}{According to \textcolor{black}{the nature of the entry values they can be found} polarizer CA~\cite{tsai2013coded}, complex-modulating CA~\cite{lai2018coded},
			phase CA~\cite{pinilla2018coded}, intensity-modulating CA~\cite{qiao2019intensity}, among others. This work \textcolor{black}{focuses on} CAs modulating the intensity through} (\romannumeral 1) binary CA (BCA)~\cite{gottesman1989new}, with opaque and translucent elements that completely block or let pass the wavefront, and (\romannumeral 2) real-valued CA~\cite{diaz2019adaptive} that attenuates the wavefront at different levels. Figure \ref{fig:proposed_method} a. shows a real-valued and a BCA whose entries represent the attenuation and block/unblock effects. The BCA provides the benefits of \textcolor{black}{lower storage requirement~\cite{duarte2008single}, higher speed during the integration time of detection~\cite{diaz2015high},} and easier fabrication and calibration, with the current technology, in comparison to the real-valued CA~\cite{mejia2018binary}. 
	 
The projected encoded measurements can be used to recover different fields depending on the modulation. To name, compressive spectral imaging systems (CSI) modulate the spatial-spectral field acquiring multiple projections of the scene~\cite{cao2016computational}; compressive imaging systems (CIS) modulate the spatial field acquiring a set of multiplexed versions to recover a gray-scale image~\cite{duarte2008single}; coded polarization imaging systems (CPI) modulate the polarization field trough a micropolarizer array~\cite{tu2017division}; coded depth imaging systems (CDI) acquire depth-dependent blurred versions of the scene~\cite{levin2007image}; coded diffraction patterns (CDP) modulate the phase and amplitude of the scene~\cite{bacca2019super}; and compressive light-field systems (CLF) modulate the spatial-angular information and acquire 2D light-field projections~\cite{marwah2013compressive}; Figure \ref{fig:proposed_method} b. lists the above mentioned optical coding systems and relates them with the corresponding application field. 
	
The performance of computational imaging (CI) tasks from projected encoded measurements such as, recovery~\cite{bacca2021compressive}, segmentation~\cite{hinojosa2018coded}, classification~\cite{bacca2020coupled}, detection~\cite{lohit2016direct} and parameter estimation~\cite{davenport2010signal}, depends on both, the computational processing method and the sensing protocol~\cite{ongie2020deep}. Hence, a line of research focuses on computational processing methods from projected encoded measurements, commonly based on optimization problems using prior knowledge through a regularizer that penalizes the main cost function. To name,~\cite{wang2015dual} assumes low total-variation to exploit the smooth transitions in the spatial domain;~\cite{correa2016multiple} assumes sparsity to exploit the sparse representation of the scene with few non-zero coefficients in a given basis; and~\cite{fu2016exploiting,wang2016adaptive,bacca2019noniterative,gelvez2017joint,gelvez2020nonlocal} assume low-rank to exploit the high structural correlations. Nonetheless, prior based methods do not represent the wide variety and non-linearity of the underlying scene, \textcolor{black}{in which the entire reconstruction} \textcolor{black}{of the image} from the projected encoded measurements \textcolor{black}{is commonly required to then achieve an acceptable quality when performing} high-level tasks. Recent data-driven deep-learning (DL) approaches employ the growing amount of available datasets to learn a non-linear transformation that maps the projected encoded measurements to the desired output, such that, high-level tasks are directly addressed from the projected encoded measurements by easily changing the deep neural network (DNN) architecture and the loss function~\cite{bacca2020coupled,liu2012texture}.

Meanwhile, another line of research focuses on the sensing protocol, particularly, on the design of the light wavefront response at each entry of the CA, known as spatial distribution. For this,~\cite{correa2015snapshot, mejia2018binary} adopt hand-designed assumptions based on prior knowledge of the sensing protocol; and~\cite{arguello2014colored, cuadros2014coded, elad2007optimized,hong2018efficient} adopt theoretical constraints as mutual coherence and concentration of measure. The design of a BCA is preferred over a real-valued CA due to its easy implementation and calibration in the physical system. Nonetheless, the mathematical design of BCAs 
\textcolor{black}{connotes the solution of an integer programming problem whose complexity is NP-complete~\cite{karp1972reducibility}, where efficient gradient-descent-based methods do not guarantee lying in a feasible solution at each iteration.} 
Hence, most of previous methods relax the complexity by applying a thresholding operation; the process of designing binary elements in the CA is commonly known as binarization. Further, most of the methods adjust empirically critical assembling properties, such as the amount of light that passes through the CA, known as transmittance\cite{marcuse1972light}, the number of captured projections with different CAs, known as snapshots, the correlation between the CAs along the snapshots, and the sensing matrix conditionality, or significant modeling factors such as, manufacturing noise of CA, and number of trainable parameters, 
which together to the spatial distribution greatly determine the performance of optical coding systems~\cite{arce2014compressive}.



Thus, at one end, DL methods employ accessible datasets to perform specific CI tasks independently of the sensing protocol determined by the CA, and at the other end, designed CAs surpass non-designed CAs, termed random CAs, in terms of storage, speed, and quality, but do not consider the available data and the imaging task. Aware of the importance of both, a novel end-to-end deep learning (E2E) approach focuses on simultaneously learning the parameters of the sensing protocol and the DNN by designed a fully-differentiable forward model of the optical coding system which is incorporated as a layer in the DNN to conduct any particular imaging task. To name,~\cite{li2020jointly,wu2019learning,mdrafi2020joint,mousavi2017deepcodec,tran2020multilinear} model the sensing matrix as a fully-connected layer, coupled to a DNN for recovery or classification tasks. Nonetheless, most of these learned sensing matrices do not fit with the structure of implementable optical systems, especially in the case of BCAs. Hence,~\cite{wang2018hyperreconnet,horisaki2020deeply,iliadis2020deepbinarymask,li2020end,shi2019image} include a piece-wise threshold function at the end of each forward step,~\cite{fu2020single} models the CA with the sigmoid activation function, and~\cite{higham2018deep, bacca2020coupled} include a regularizer to address the binarization and obtain implementable sensing matrices. Although the latter strategy presents a more versatile approach to address the binarization, it does not consider the whole important assembling properties and modeling considerations, so that, the search of functions to regularize CA designs is an open field. 

A reliable E2E approach then requires the selection of four fronts: a suitable optical coding system, the correct modeling of the CA, the accomplishment of a CI task, and the training of a DNN with proper loss and regularizer functions. However, previous methods present some gaps. \textcolor{black}{First, they do not allow to deal with all four fronts at the same time and they can not be easily adapted for different purposes since they have been developed for a specific application at the time.} Second, they do not consider the entire assembling properties of CAs, as transmittance, number of snapshots, correlation between snapshots, and sensing matrix conditionality. Third, they do not model issues as the number of trainable parameters and manufacturing noise. Therefore, we propose a customizable E2E approach
that jointly designs the CA and learns a DNN for CI tasks which intrinsically employs specific fields of the plenoptic function as shown in Fig.~\ref{fig:proposed_method}~c. This general scheme is based on the incorporation of strategic regularizers in the loss function that accomplish particular properties of the sensing protocol design, especially of the CA. Then, we formulate a set of functions to decide the transmittance level, number of snapshots, correlation, and sensing matrix conditionality, and a family of functions to generate binary and real-valued CAs. We further introduce a dynamic strategy to choose the associated regularization parameter. Finally, we analyze the manufacturing noise and the number of trainable parameters in real setups.
We validate the effectiveness of the proposal along three optical coding systems; three categories of images: gray-scale images, spectral images, and 3D images; and three CI tasks: restoration, classification and semantic segmentation. 

	\section{Mathematical Considerations}
	
	
The plenoptic function represents the electromagnetic wavefront emitted by an object when illuminated with a light source~\cite{mcmillan1995plenoptic}. Assuming that a ray carries optical energy, this function can be modeled with at least eight dimensions as $f(x,y,z,\theta,\psi,t,\lambda,p)$, where $(x,y)$ stand for spatial, $z$ for depth, $(\theta,\psi)$ for angular views, $t$ for temporal, $\lambda$ for wavelength, and $p$ for polarization dimension. The acquisition of all dimensions is challenging because current optical systems rely on converting photons to electrons, i.e., they measure the photon flux per unit of surface area. Then, CAs are incorporated to encode a subset of the dimensions of the plenoptic function before being projected to the sensor.

	\subsection{Coded Aperture Ensembles}
	\label{sec:coded_model}
Considering the plenoptic function, the sensing protocol and CA can be modeled as 8-dimensional operators. This work focuses on protocols that modulate 2 or 3 dimensions at the time, where the CAs are modeled as 2D or 3D arrays. These models assume that the light response remains approximately constant over a spatial square region of size $\Delta_p \times \Delta_p$, named pixel, having a specific response to the incoming light. 
	
	
	A 2D CA, $\textcolor{black}{\boldsymbol{\hat{\Phi}}(x,y)} \in \mathbb{R}^{M \times N}$, mostly models the modulation of the spatial field of the plenoptic function determined by the following transmittance function 
	\begin{equation}
	\textcolor{black}{\boldsymbol{\hat{\Phi}}(x,y)} = \sum_{i,j} \boldsymbol{\Phi}_{i,j} \text{rect} \left(\frac{x}{\Delta_p} -i ,\frac{y}{\Delta_p} -j \right),
	\label{eq:2D_reoresentation}
	\end{equation}
	for $i=1,\ldots, M$, $j=1, \ldots, N$ spatial pixels, where $\text{rect}(\cdot)$ denotes the rectangular step function. When the continuous light response at the $(i,j)^{th}$ discrete location is modeled with the intention of generating a physical binary behavior, i.e. $\boldsymbol{\Phi}_{i,j} \in \{0,1\}$ to block ($0$) or unblock ($1$) the entire composition of the light, this model is referred to as BCA. Unlike, when it is modeled as $\boldsymbol{\Phi}_{i,j}\in[0,1]$, this model is referred to as real-valued CA.
	
\begin{figure}[b!]
		 	\begin{tikzpicture}
		 		 \tikzstyle{every node}=[font=\small]
		 		\node[inner sep=0pt] (Examples) at (0\linewidth,0)
 {\includegraphics[width=0.98\linewidth]{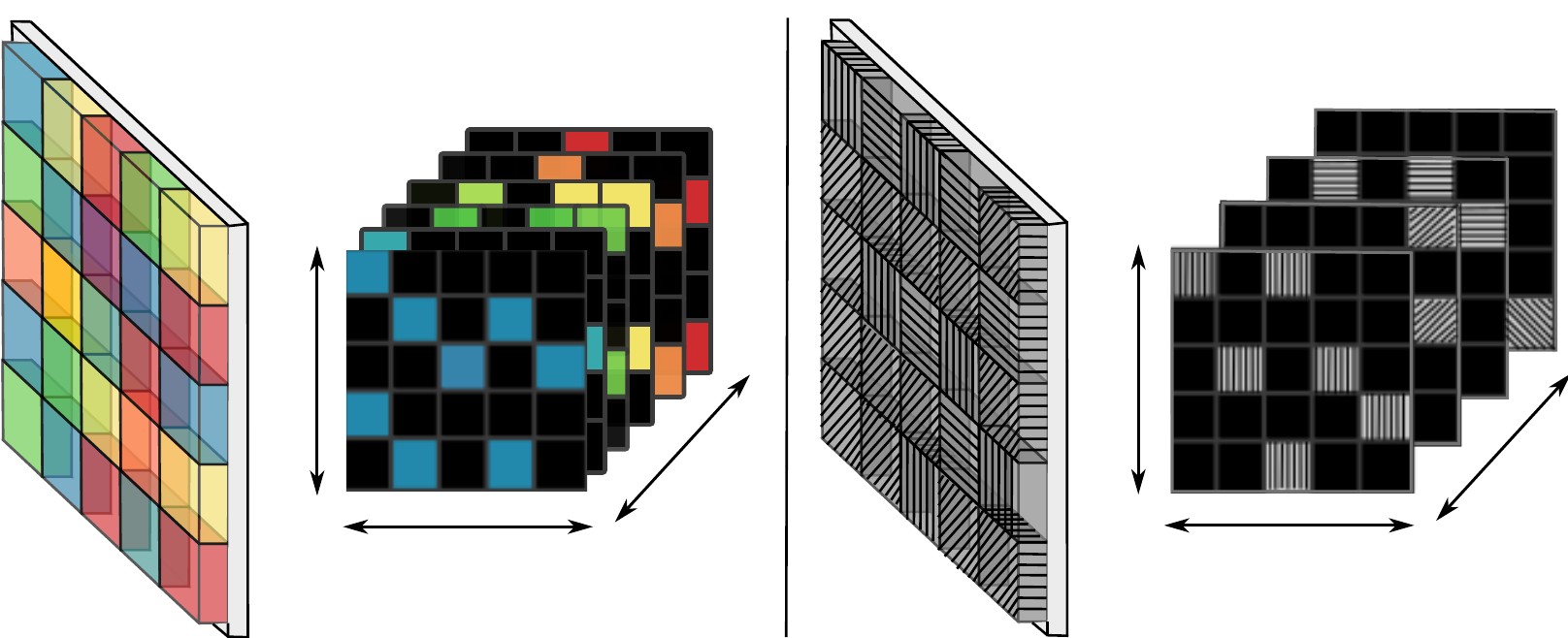}};
\node at (-0.25\linewidth,1.5cm) {\textbf{a. Colored CA}};
\node at (-0.32\linewidth,-0.3cm) {$x$};
\node at (-0.2\linewidth,-1.4cm) {$y$};
\node at (-0.04\linewidth,-0.75cm) {$\lambda$};
\node at (-0.43\linewidth,-2.0cm) {Physical};
\node at (-0.43\linewidth,-2.3cm) {scheme};
\node at (-0.17\linewidth,-2.0cm) {Mathematical};
\node at (-0.17\linewidth,-2.3cm) {model};

\node at (0.25\linewidth,1.5cm) {\textbf{b. Polarizer CA}};
\node at (0.20\linewidth,-0.3cm) {$x$};
\node at (0.32\linewidth,-1.4cm) {$y$};
\node at (0.48\linewidth,-0.75cm) {$p$};

\node at (0.09\linewidth,-2.0cm) {Physical};
\node at (0.09\linewidth,-2.3cm) {scheme};
\node at (0.35\linewidth,-2.0cm) {Mathematical};
\node at (0.35\linewidth,-2.3cm) {model};	\end{tikzpicture}\vspace{-1.1em}
			\caption{Visual representation of two 3D CA. (a.) A colored CA whose third dimension discretizes the wavelength. (b.) A polarizer CA whose third dimension discretizes the polarization angles.}
		\label{fig:CAs}
	\end{figure}
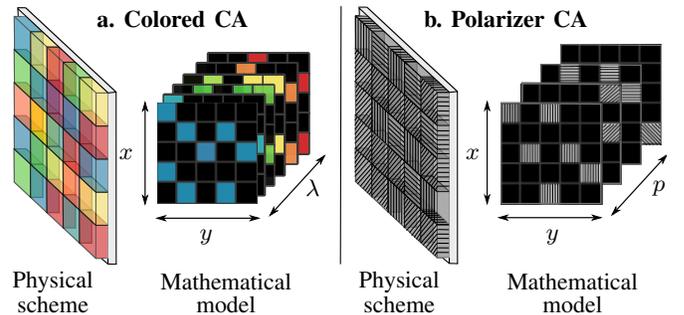

A 3D CA, $\textcolor{black}{\boldsymbol{\hat{\Phi}}(x,y,\zeta)} \in \mathbb{R}^{M \times N \times L}$, mostly models the 2D spatial distribution along a third dimension of the plenoptic function, determined by the manufacturing, configuration, or distance where the CA is placed with respect to the sensor. To name, the colored-CA models the 2D spatial distribution along the 1D spectrum, i.e. it contains the spectral response of various optical filters at each spatial location of the CA~\cite{arguello2014colored}; the polarizer-CA, models the 2D spatial distribution along the 1D polarization angle, i.e. it is a micro-polarizer array composed of small apertures of wire grid at different polarization angles~\cite{nordin1999micropolarizer}; the temporal-CA models the 2D spatial distribution along the time, i.e. it is a 2D CA whose spatial distribution changes faster against the integration time of the sensor~\cite{baraniuk2017compressive}; and the depth-CA models the 2D spatial distribution along the depth, i.e it is a set of 2D CAs placed between the lenses and the sensor to encode the depth and angular views~\cite{marwah2013compressive}. Mathematically, a 3D CA is given by
	\begin{equation}
	\textcolor{black}{\boldsymbol{\hat{\Phi}}(x,y,\zeta)} = \sum_{i,j,\ell} \boldsymbol{\Phi}_{i,j,\ell} \text{rect} \left(\dfrac{x}{\Delta_p} -i ,\dfrac{y}{\Delta_p} -j, \dfrac{\zeta}{\Delta_\zeta}-\ell \right),
	\label{eq:Coded}
	\end{equation}
 where $\zeta$ denotes the plenoptic function variable of interest, with $L$ discretizations of size $\Delta_\zeta$, and $\boldsymbol{\Phi}_{i,j,\zeta}$ is the wavefront response which can be binary $\boldsymbol{\Phi}_{i,j,\zeta } \in \{ 0,1\}$, or real-valued $\boldsymbol{\Phi}_{i,j,\zeta } \in [0,1]$, depending on the used technology. Figure~\ref{fig:CAs} shows two 3D CA distinguishing between the physical scheme and the 3D mathematical structure that models the CA.
 
	\subsection{Coded Sensing Observation Model}
	
	The coded sensing observation model refers to the forward model to acquire a set of projected encoded measurements through any optical coding system independently of the nature and assembling properties of the CA. 
	Mathematically, the projected encoded measurements obtained in one snapshot can be expressed as
	\begin{equation}
	\mathbf{g}=\mathbf{H}_{\boldsymbol{\Phi}}\mathbf{f} + \boldsymbol{\eta},
	\label{eq:sensing_model}
	\end{equation}
	where $\mathbf{g}\in \mathbb{R}^{m}$ denotes the projected encoded measurements, $\mathbf{f}\in \mathbb{R}^{n}$ denotes the underlying scene, $\mathbf{H}_{\boldsymbol{\Phi}}\in\mathbb{R}^{m\times n}$ models the sensing matrix whose structure is determined by the setup of the optical coding system and the corresponding CA ($\boldsymbol{\Phi}$), and $\boldsymbol{\eta}\in \mathbb{R}^{m}$ stands for the noise. To highlight, the CA is the main customizable physical element in an optical coding system, so that, it highly determines the performance of the CI task when using projected encoded measurements. Further, some optical coding systems enable the acquisition of multiple snapshots of the same scene, assuming that the scene remains constant over a time-lapse, by easily varying the used CA. This process is referred to as multishot acquisition~\cite{duarte2008single}.

	Mathematically, for a number of $S$ snapshots, each projected encoded measurement $\{\mathbf{g}^{s}\}_{s=1}^{S}$ is obtained with a different sensing matrix $\{\mathbf{H}_{\boldsymbol{\Phi}^{s}}\}_{s=1}^{S}$, modeled as in \eqref{eq:sensing_model}. The multishot acquisition process can be compactly expressed as
	\begin{equation}
	 \tilde{\mathbf{g}}=\tilde{\mathbf{H}}_{\tilde{\boldsymbol{\Phi}}}\mathbf{f} + \tilde{\boldsymbol{\eta}},
	 \label{eq:all_measurments}
	\end{equation}
	where $\tilde{\mathbf{g}}=[(\mathbf{g}^1)^T, \cdots, (\mathbf{g}^{S})^T]^T$ stacks the projected encoded measurements, $\tilde{\mathbf{H}}_{\tilde{\boldsymbol{\Phi}}} = [(\mathbf{H}_{\boldsymbol{\Phi}^{1}})^T, \cdots, (\mathbf{H}_{\boldsymbol{\Phi}^{S}})^T]^T$ stacks the corresponding sensing matrices, $\tilde{\boldsymbol{\Phi}} = \{\boldsymbol{\Phi}^s\}_{s=1}^{S}$ is a set containing the corresponding CAs, and $\boldsymbol{\eta}=[(\tilde{\boldsymbol{\eta}}^1)^T, \cdots, (\tilde{\boldsymbol{\eta}}^{S})^T]^T$ stacks the measurements noise.
	The ratio between the amount of observed measurements and the size of the underlying scene is known as compression ratio given by $\gamma = Sm/n$.
	
	\section{Proposed End-To-End approach}
	\label{sec:Design}
	Our proposal aims to couple the design of the sensing matrix $\tilde{\mathbf{H}}_{\boldsymbol{\tilde{\Phi}}}$ together with the achievement of a CI task of interest using an E2E approach. Thus, it jointly optimizes the set of CAs $\tilde{\boldsymbol{\Phi}}$, and the parameters $\boldsymbol{\theta}$, of a chosen DNN $\mathcal{M}_{\boldsymbol{\theta}}(\cdot)$, by minimizing a cost function composed by the sum of the loss function $ \mathcal{L}_
	{task}(\cdot,\cdot)$, to achieve the task, and a customizable regularization function $R_q(\tilde{\boldsymbol{\Phi}})$, that promotes particular properties in the set of CAs. Mathematically, 
	given a set of $K$ scenes $\{\mathbf{f}_{k}\}_{k = 1}^{K}$, and their corresponding task outputs $\{\mathbf{d}_k\}_{k = 1}^{K}$, our coupled optimization problem is formulated as 
		\begin{align}
	\{\tilde{\boldsymbol{\Phi}}^*,\boldsymbol{\theta}^*\} \hspace{2mm} \in & \hspace{2mm} \argmin_{\tilde{\boldsymbol{\Phi}}, \boldsymbol{\theta}} \hspace{2mm} \mathcal{L}(\tilde{\boldsymbol{\Phi}}, \boldsymbol{\theta} | \mathbf{f}, \mathbf{d}), \label{eq:regularization1} \\ 
	\nonumber \mathcal{L}(\tilde{\boldsymbol{\Phi}}, \boldsymbol{\theta} | \mathbf{f}, \mathbf{d}) = &  \frac{1}{K} \sum_{k} \mathcal{L}_
	{task}\left(\mathcal{M}_{\boldsymbol{\theta}}(\tilde{\mathbf{H}}_{\tilde{\boldsymbol{\Phi}}}\mathbf{f}_k ), \mathbf{d}_k\right)+ \rho_q R_q(\tilde{\boldsymbol{\Phi}}),
	\end{align}
 where, $\rho_q > 0$ is a regularization parameter. 

Then, we follow a two-module procedure that models the encoder and decoder steps to solve \eqref{eq:regularization1}. The first, referred to as the sensing protocol module, consists of an optical layer that learns the sensing matrix, $\tilde{\mathbf{H}}_{\tilde{\boldsymbol{\Phi}}}$, with the set of CAs $\tilde{\boldsymbol{\Phi}}$, as the trainable parameters. Notice that $\tilde{\boldsymbol{\Phi}}$ is regularized by the function $R_q(\tilde{\boldsymbol{\Phi}})$, which aims to guarantee the formation of CAs that satisfy specific implementability, transmittance, number of snapshots, correlation between snapshots, and conditionality constraints. The optical layer is directly connected to the second module, referred to as the task module that consists of a chosen trainable DNN with various hidden layers to conclude a task. Figure \ref{fig:proposed_method} outlines the coupled E2E proposed approach. It can be seen that the CAs directly affect the task and vice-versa, since the backward step to update the values of $\tilde{\boldsymbol{\Phi}}$ takes into account the trade-off between the error given by the loss of the task, and the regularizer function of the CAs.

Remark that once the set of CAs is optimized for the task, the designed sensing matrix $\tilde{\mathbf{H}}_{\tilde{\boldsymbol{\Phi}}}^*$ can be used to acquire new projected encoded measurements ($\bar{\mathbf{g}}$), and the pre-trained DNN be applied to estimate the task output $\bar{\mathbf{d}}$, as	$\bar{\mathbf{d}} = \mathcal{M}_{\boldsymbol{\theta}^*}(\bar{\mathbf{g}}).$ Thus, the pretrained DNN is used as an inference operator.
	
	Next subsections detail the constraints for designing CAs of a coding sensing protocol with the proposed regularizers. 

	\subsection{Binary Coded Aperture Implementation Constraint}
	\label{sub:BCA}
\begin{figure}[t!]
	\begin{tikzpicture}
	 \tikzstyle{every node}=[font=\small]
	\node[inner sep=0pt] (Family) at (0\linewidth,0){
		\includegraphics[width=1\linewidth]{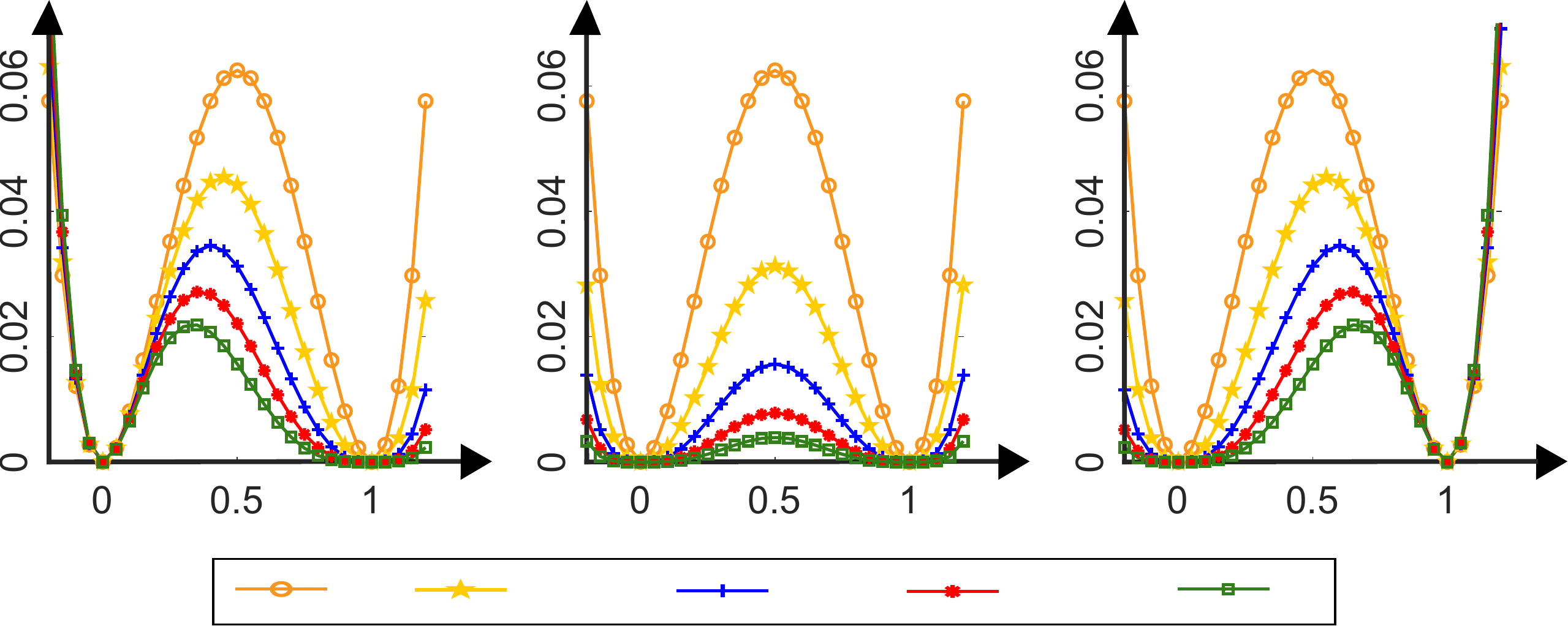}};
		
\node at (-0.44\linewidth,1.8cm) {\footnotesize{$R(\boldsymbol{\Phi}_{i,j})$}};
\node at (-0.1\linewidth,1.8cm) {\footnotesize{$R(\boldsymbol{\Phi}_{i,j})$}};
\node at (0.24\linewidth,1.8cm) {\footnotesize{$R(\boldsymbol{\Phi}_{i,j})$}};

\node at (-0.33\linewidth,2.2cm) {{ \textbf{a.} ${p_2}>{p_1} = 1 $}};
\node at (0.0\linewidth,2.2cm) {{\textbf{b.} ${p_2}={p_1}$}};
\node at (0.33\linewidth,2.2cm) {{\textbf{c.} ${p_2}=1<{p_1}$}};

\node at (-0.27\linewidth,-1.6cm) {$1$};
\node at (-0.13\linewidth,-1.6cm) {$1.25$};
\node at (0.02\linewidth,-1.6cm) {$1.5$};
\node at (0.18\linewidth,-1.6cm) {$1.75$};
\node at (0.33\linewidth,-1.6cm) {$2$};

\node at (-0.2\linewidth,-1.2cm) {\footnotesize{$\boldsymbol{\Phi}_{i,j}$}};
\node at (0.14\linewidth,-1.2cm) {\footnotesize{$\boldsymbol{\Phi}_{i,j}$}};
\node at (0.48\linewidth,-1.2cm) {\footnotesize{$\boldsymbol{\Phi}_{i,j}$}};
		\end{tikzpicture}\vspace{-1em}
		\caption{Family of functions to regularize BCA along different values of the tuple ${p_1},{p_2}$. The legends refer to the value of ${p_2}$ in (\textbf{a.}), ${p_1}={p_2}$ in \textbf{(b.)}, and ${p_1}$ in \textbf{(c.)}.}
		\label{fig:family}
	\end{figure}
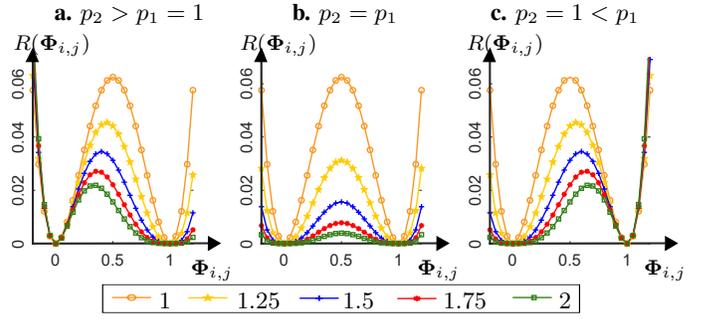

	The binarization constraint is addressed by proposing a family of functions whose minima are obtained uniquely when the elements of the CAs are either ($0$) or ($1$). Then, it is incorporated in \eqref{eq:regularization1} as the regularizer $ R_1(\tilde{\boldsymbol{\Phi}})$ given by
	\begin{equation}
	 R_1(\tilde{\boldsymbol{\Phi}}) = \frac{1}{S} \sum_s \sum_{i,j,{\ell}} \left(({\boldsymbol{\Phi}_{i,j,{\ell}}^s})^2\right)^{p_1}\left(({\boldsymbol{\Phi}_{i,j,{\ell}}^s}-1)^2\right)^{p_2},
	\label{eq:family_binary_0}
	\end{equation}
	where $p_1,p_2 \in \mathbb{R_{++}}$ are two hyper-parameters that provide variability in the function curve as illustrated in Fig. \ref{fig:family}, where
 the behaviour of the family of functions along three different combinations of the tuple $({p_1},{p_2})$ is presented. In particular, if ${p_1}={p_2}$ the graph is symmetric, which turns in equal speed to converge to any of the minima when minimizing the function. Meanwhile, when ${p_2}>{p_1}$ or ${p_1}>{p_2}$, the function presents a bias that leads to a faster direction to converge to one of the minima, $0$ or $1$. Introducing ${p_1}$ and ${p_2}$ to control the bias will be useful to determine the transmittance of the designed CAs.

	Further, works in~\cite{wagadarikar2008single,candes2008introduction,higham2018deep} address the binarization constraint through CAs composed by $\pm 1$ values\footnote{Codification with negative values can be physically achieved by estimating and subtracting the mean light intensity from each measurement, which can be achieved with a full-one CA~\cite{duarte2008single} \textcolor{black}{For more physical details see~\cite{bacca2020coupled}}.}. We address this constraint by proposing a family of functions whose minimums are found uniquely when the elements of the CAs are either $-1$ and $1$. Hence, the proposed family of functions is incorporated in \eqref{eq:regularization1} as the regularizer $ R_2(\tilde{\boldsymbol{\Phi}})$ given by
	\begin{equation}
	 	 R_2(\tilde{\boldsymbol{\Phi}}) = \frac{1}{S} \sum_s \sum_{i,j,{\ell}} \left((\boldsymbol{\Phi}_{i,j,{\ell}}^{{s}}+1)^2\right)^{{p_1}}\left((\boldsymbol{\Phi}_{i,j,{\ell}}^{{s}}-1)^2\right)^{{p_2}},
	\label{eq:binary_1}
	\end{equation}
	which generalizes the function presented in~\cite{bacca2020coupled,higham2018deep}. Particularly, when ${p_1}={p_2}=1$, \eqref{eq:binary_1} produces the function used in~\cite{higham2018deep} for reconstruction and in~\cite{bacca2020coupled} for classification. 

	\subsection{Real-valued Coded Aperture Implementation Constraint}
	\label{subsec:gray-scale}
	The range of admitted values in a real-valued CA ($\boldsymbol{\Phi}\in [0,1]$) is wider than in the binary one, which mathematically results in a more relaxed constraint. However, the wider the range of values in the mathematical design, the harder the fabrication and calibration in the physical system~\cite{diaz2019adaptive}. This difficulty can be alleviated by fixing a set $\{\kappa_d\}_{d=1}^{D} \in \mathbb{R}$ of different quantization levels to appear in the design. We address this real-valued constraint by proposing a family of functions whose minima are obtained uniquely when the elements of the CAs are the fixed quantization levels $\kappa_d$. Then, we incorporate it in \eqref{eq:regularization1} as the regularizer $ R_3(\tilde{\boldsymbol{\Phi}})$ given by
	\begin{equation}
	 { R_3(\tilde{\boldsymbol{\Phi}}) = \frac{1}{S} \sum_s} \sum_{i,j,{\ell}} \prod_{d} \left(\left( \boldsymbol{\Phi}_{i,j,{\ell}}^{{s}} - \kappa_d \right)^{2}\right)^{p_d},
	\label{eq:family_gray_0}
	\end{equation}
	where $\{p_d\}_{d=1}^{D} \in \mathbb{R_{++}}$ are the hyper-parameters associated for each quantization level. Notice that \eqref{eq:family_gray_0} is a polynomial of degree $\prod_{d}{2p_d}$ whose range is always positive and whose roots correspond to the targeted quantization levels in the CA design. Thus, \eqref{eq:family_gray_0} is the generalization of our previous family of functions in \eqref{eq:family_binary_0} and \eqref{eq:binary_1} which are the particular cases of 2 target values, i.e., $D=2$ with $\kappa_1 = 0, \kappa_2 = 1$ for \eqref{eq:family_binary_0}, and $\kappa_1 = -1, \kappa_2 = 1$, for \eqref{eq:binary_1}, respectively.

	\subsection{Coded Aperture Transmittance Constraint}
	\label{sub:trans}
	The transmittance level is a crucial property that affects the calibration of the optical coding system and determines the proper utilization of the light in the acquired measurements. Therefore, adjusting the transmittance level is an important step to accomplish a task. For instance, in spectral imaging, a high-transmittance level is desired to increase the signal to noise power ratio~\cite{rueda2015multi}, unlike, in X-ray tomography, a low-transmittance level is desired to minimize radiation to the objects~\cite{mojica2017high}. Specifically, varying the transmittance level unbalances the distribution of the quantization levels in the CA and produces an ill-conditioned sensing matrix; in consequence, the performance of the system can be affected. The general family of functions can indirectly control the transmittance level adjusting the hyperparameter $\kappa_d$. Nonetheless, to achieve an exact targeted value we address this transmittance level constraint by proposing the following regularization function $ R_4(\tilde{\boldsymbol{\Phi}})$ that adjusts the transmittance level while affecting the less the possible the performance of the system 
	\begin{equation}
 { R_4(\tilde{\boldsymbol{\Phi}}) = \frac{1}{S} \sum_s} \left(\frac{\sum_{i,j,{\ell}}\boldsymbol{\Phi}_{i,j,{\ell}}^{{s}}}{MN{L}} - T_r\right)^2,
\label{eq:transmittance_fixed}
	\end{equation}
	where $T_r \in [0,1]$ is a customizable hyperparameter that denotes the targeted transmittance level, with $0$ and $1$ indicating to block or unblock all the incoming light.

	\subsection{Number of Snapshots Constraint}
		\label{sub:numberSnapshots}
The aim of acquiring multiple snapshots is to efficiently increase the amount of observed information related to the properties of the scene improving the task performance. However, the number of snapshots implies a trade-off between the task performance and the needed time for acquiring and processing more snapshots~\cite{baraniuk2017compressive}. Therefore, it is essential to determine the least amount of optimal snapshots $S$ that achieve the highest performance. We propose to address this number of snapshots constraint by using the following regularizer
	\begin{equation}
	R_5(\tilde{\boldsymbol{\Phi}}) := \sum_{s'}^{S'} \sqrt{\sum_{i,j,{\ell}} \left(\boldsymbol{\Phi}_{i,j,{\ell}}^{{s'}}\right)^2},
	\label{eq:regu_shots} 
	\end{equation}
	where $S' \geq S$ is an upper bound for the number of snapshots determined by the user. This equation can be seen as the $\ell_1$-norm through the snapshot dimension of the $\ell_2$-norm of the vectorization of each CA. This formulation is based on the traditional $\ell_{2,1}$-norm applied on matrices, which has been demonstrated to encourage all values in a column to be zero~\cite{eldar2009block,lohit2019rank}. Thus, \eqref{eq:regu_shots}, promotes all entries in a CA to be equal to zero, i.e., implying the no acquisition of the snapshots. 

	Notice that expression in \eqref{eq:regu_shots} is not differentiable when all values of a CA are zero ~\cite{lohit2019rank}. Hence, we employ the sub-gradient below in the backward step 
	\begin{equation}
	\label{eq:_gradient}
	\frac{\partial R_5(\tilde{\boldsymbol{\Phi}})}{\partial\boldsymbol{\Phi}_{i,j,{\ell}}^{s'}} = \hspace{-0.0em}\left\{ \hspace{-0.2em} \begin{array}{lcc}
	\frac{\boldsymbol{\Phi}_{i,j,{\ell}}^{s'}}{\varphi(\boldsymbol{\Phi}_{i,j,{\ell}}^{s'})}, & \varphi(\boldsymbol{\Phi}_{i,j,{\ell}}^{s'}) \neq 0 
	\\ 0, & \mbox{otherwise},
	\end{array}
	\right. 
	\end{equation}
	where $\varphi(\boldsymbol{\Phi}_{i,j,{\ell}}^{s'}) = \sqrt{\sum_{i,j,{\ell}}
	(\boldsymbol{\Phi}_{i,j,{\ell}}^{s'})^2}$.

	Notice that regularizers proposed up to this section are directly related to assembling properties of the CA, such as the implementability of the obtained quantization levels, the adjustment of the transmittance, and the selection of the number of snapshots to be acquired. Unlike, the following two regularizers are proposed in order to achieve a better performance in the solution of the CI task.
		
	\subsection{Multishot Coded Aperture Correlation Constraint}
	\label{sub:correlation}
	The correlation between the CAs $\{\boldsymbol{\Phi}^s\}_{s=1}^{S}$ used in a multishot acquisition scheme is crucial to increase the observed information of the underlying scene. Specifically, it has been demonstrated in~\cite{correa2015snapshot} that the less correlated the CAs, the greater the amount of acquired with less snapshots. We address the correlation constraint by using a generalized function that minimizes the correlation between the $S$ designed CAs. Then, we incorporate it in \eqref{eq:regularization1} as the regularizer $ R_6(\tilde{\boldsymbol{\Phi}})$ given by
	
	\begin{equation}
	R_6(\tilde{\boldsymbol{\Phi}}) =\frac{\sum_{i,j,{\ell}} \left(\prod_{s} \boldsymbol{\Phi}_{i,j,\ell}^{s} \right)}{MNL}.
	\label{eq:correlation}
	\end{equation}
	Observe that for $S=2$, expression in \eqref{eq:correlation} results in the numerator part of the Pearson correlation for two CAs~\cite{benesty2009pearson}.


\begin{figure*}[!t]
		\begin{tikzpicture}
		 \tikzstyle{every node}=[font=\small]
		 \filldraw [fill=blue!4, draw=black , rounded corners, thick] (0,0) rectangle (1\linewidth,-3.3); 
\filldraw [fill=white, draw=black,, rounded corners] (0.01\linewidth,-0.5) rectangle (0.05\linewidth,-0.8);
\filldraw [fill=white, draw=black, rounded corners] (0.01\linewidth,-0.8) rectangle (0.05\linewidth,-1.1);
\filldraw [fill=white, draw=black,rounded corners] (0.01\linewidth,-1.1) rectangle (0.05\linewidth,-1.4);
\filldraw [fill=white, draw=black, rounded corners] (0.01\linewidth,-1.4) rectangle (0.05\linewidth,-1.7);
\filldraw [fill=white, draw=black, rounded corners] (0.01\linewidth,-1.7) rectangle (0.05\linewidth,-2.0);
\filldraw [fill=white, draw=black, rounded corners] (0.01\linewidth,-2.0) rectangle (0.05\linewidth,-2.3);
\filldraw [fill=white, draw=black, rounded corners] (0.01\linewidth,-2.3) rectangle (0.05\linewidth,-2.6);
\filldraw [fill=white, draw=black, rounded corners] (0.01\linewidth,-2.6) rectangle (0.05\linewidth,-2.9);
\filldraw [fill=white, draw=black, rounded corners] (0.01\linewidth,-2.9) rectangle (0.05\linewidth,-3.2);

\node at (0.03\linewidth,-3.05){G};
\node at (0.03\linewidth,-2.75){E};
\node at (0.03\linewidth,-2.45){D};
\node at (0.03\linewidth,-2.15){C};
\node at (0.03\linewidth,-1.85){F};
\node at (0.03\linewidth,-1.55){H};
\node at (0.03\linewidth,-1.25){I};
\node at (0.03\linewidth,-0.95){B};
\node at (0.03\linewidth,-0.65){A};
\node at (0.03\linewidth,-0.25) {\textbf{Exp.}};

\filldraw [fill=white, draw=black, rounded corners] (0.13\linewidth,-0.5) rectangle (0.28\linewidth,-1.25);
\filldraw [fill=white, draw=black, rounded corners] (0.13\linewidth,-1.25) rectangle (0.28\linewidth,-2.0);
\filldraw [fill=white, draw=black,rounded corners] (0.13\linewidth,-2.0) rectangle (0.28\linewidth,-2.75);

\node at (0.13\linewidth,-2.55) [anchor=west] {regularizer in \eqref{eq:family_gray_0}};
\node at (0.13\linewidth,-2.25) [anchor=west] {real-valued $[0,1]$};
\node at (0.13\linewidth,-1.8) [anchor=west] {regularizer in \eqref{eq:binary_1}};
\node at (0.13\linewidth,-1.5) [anchor=west] {Binary $\{-1,1\}$};
\node at (0.13\linewidth,-1.05) [anchor=west] {regularizer in \eqref{eq:family_binary_0}};
\node at (0.13\linewidth,-0.75) [anchor=west] {Binary $\{0,1\}$ with};
\node at (0.205\linewidth,-0.25) {\textbf{1. Coded aperture}};

\filldraw [fill=white, draw=black, rounded corners] (0.36\linewidth,-0.5) rectangle (0.47\linewidth,-1.25);
\filldraw [fill=white, draw=black, rounded corners] (0.36\linewidth,-1.25) rectangle (0.47\linewidth,-2);
\filldraw [fill=white, draw=black,rounded corners] (0.36\linewidth,-2) rectangle (0.47\linewidth,-2.75);

\node at (0.36\linewidth,-2.375) [anchor=west] {C-Depth~\cite{levin2007image}};
\node at (0.36\linewidth,-1.625) [anchor=west] {SPC~\cite{duarte2008single}};
\node at (0.36\linewidth,-0.875) [anchor=west] {CASSI~\cite{wagadarikar2008single}};
\node at (0.415\linewidth,-0.25) {\textbf{2. Coding system}};

\filldraw [fill=white, draw=black, rounded corners] (0.55\linewidth,-0.5) rectangle (0.67\linewidth,-1.25);
\filldraw [fill=white, draw=black, rounded corners] (0.55\linewidth,-1.25) rectangle (0.67\linewidth,-2);
\filldraw [fill=white, draw=black,rounded corners] (0.55\linewidth,-2) rectangle (0.67\linewidth,-2.75);

\node at (0.55\linewidth,-2.55) [anchor=west] {segmentation};
\node at (0.55\linewidth,-2.2) [anchor=west] {Semantic};
\node at (0.55\linewidth,-1.625) [anchor=west] {Classification};
\node at (0.55\linewidth,-0.875) [anchor=west] {Reconstruction};
\node at (0.61\linewidth,-0.25) {\textbf{3. Imaging task}};

\filldraw [fill=white, draw=black,, rounded corners] (0.75\linewidth,-0.5) rectangle (0.99\linewidth,-0.9);
\filldraw [fill=white, draw=black, rounded corners] (0.75\linewidth,-0.9) rectangle (0.99\linewidth,-1.3);
\filldraw [fill=white, draw=black,rounded corners] (0.75\linewidth,-1.3) rectangle (0.99\linewidth,-1.7);
\filldraw [fill=white, draw=black, rounded corners] (0.75\linewidth,-1.7) rectangle (0.99\linewidth,-2.1);
\filldraw [fill=white, draw=black, rounded corners] (0.75\linewidth,-2.1) rectangle (0.99\linewidth,-2.5);
\filldraw [fill=white, draw=black, rounded corners] (0.75\linewidth,-2.5) rectangle (0.99\linewidth,-2.9);

\node at (0.75\linewidth,-2.7)[anchor=west]{Transmittance level in \hspace{0.42em} \eqref{eq:transmittance_fixed}};
\node at (0.75\linewidth,-2.3)[anchor=west]{Number of snapshots in \eqref{eq:regu_shots}};
\node at (0.75\linewidth,-1.9)[anchor=west]{Trainable parameters in \hspace{-0.18em} \eqref{eq:Kernel_trainable}};
\node at (0.75\linewidth,-1.5)[anchor=west]{manufacturing noise in \hspace{0.05em} \eqref{eq:ManNois}};
\node at (0.75\linewidth,-1.1)[anchor=west]{Conditionality in \hspace{3em}\eqref{eq:conditionality}};
\node at (0.75\linewidth,-0.7)[anchor=west]{Multishot correlation in \hspace{-0.25em} \eqref{eq:correlation}};
\node at (0.87\linewidth,-0.25) {\textbf{4. Regularizer}};

\draw[brown, loosely dashdotted, thick] (0.55\linewidth,-0.61) -- (0.47\linewidth,-0.61);
\draw[brown,loosely dashdotted, thick] (0.36\linewidth,-0.61) -- (0.28\linewidth,-0.61);
\draw[brown, loosely dashdotted, thick] (0.13\linewidth,-0.61) -- (0.05\linewidth,-0.95);

\draw[brown, loosely dashdotted, thick] (0.13\linewidth,-0.61) -- (0.05\linewidth,-0.65);

\draw[green!80!black!60,thick] (0.55\linewidth,-1.16) -- (0.47\linewidth,-1.63);
\draw[green!80!black!60,thick] (0.36\linewidth,-1.63) -- (0.28\linewidth,-1.63);
\draw[green!80!black!60,thick] (0.13\linewidth,-1.63) -- (0.05\linewidth,-2.15);

\draw[orange,densely dashdotted,thick] (0.55\linewidth,-2.5) -- (0.47\linewidth,-2.5);
\draw[orange,densely dashdotted,thick] (0.36\linewidth,-2.5) -- (0.28\linewidth,-2.5);
\draw[orange,densely dashdotted,thick] (0.36\linewidth,-2.5) -- (0.28\linewidth,-1.05);
\draw[orange,densely dashdotted,thick] (0.13\linewidth,-2.25) -- (0.05\linewidth,-2.45);
\draw[orange,densely dashdotted,thick] (0.13\linewidth,-1.05) -- (0.05\linewidth,-2.45);

\draw[blue, dashed, thick] (0.75\linewidth,-2.7) -- (0.67\linewidth,-1.82);
\draw[blue, dashed, thick] (0.55\linewidth,-1.82) -- (0.47\linewidth,-1.82);
\draw[blue, dashed, thick] (0.36\linewidth,-1.82) -- (0.28\linewidth,-1.82);
\draw[blue, dashed, thick] (0.13\linewidth,-1.82) -- (0.05\linewidth,-2.75);

\draw[violet, dotted, thick] (0.75\linewidth,-2.3) -- (0.67\linewidth,-1.05);
\draw[violet, dotted, thick] (0.55\linewidth,-1.05) -- (0.47\linewidth,-1.44);
\draw[violet, dotted, thick] (0.36\linewidth,-1.44) -- (0.28\linewidth,-0.94);
\draw[violet, dotted, thick] (0.13\linewidth,-0.94) -- (0.05\linewidth,-1.85);

\draw[red, loosely dotted, ultra thick] (0.75\linewidth,-1.9) -- (0.67\linewidth,-0.94);
\draw[red, loosely dotted, ultra thick] (0.55\linewidth,-0.94) -- (0.47\linewidth,-0.94);
\draw[red, loosely dotted, ultra thick] (0.36\linewidth,-0.94) -- (0.28\linewidth,-2.25);
\draw[red, loosely dotted, ultra thick] (0.13\linewidth,-2.5) -- (0.05\linewidth,-3.05);

\draw[cyan, dashdotted, thick] (0.75\linewidth,-1.5) -- (0.67\linewidth,-0.83);
\draw[cyan, dashdotted, thick] (0.55\linewidth,-0.83) -- (0.47\linewidth,-0.83);
\draw[cyan, dashdotted, thick] (0.36\linewidth,-0.83) -- (0.28\linewidth,-0.83);
\draw[cyan, dashdotted, thick] (0.13\linewidth,-0.83) -- (0.05\linewidth,-1.55);

\draw[magenta, densely dashed, thick] (0.75\linewidth,-0.7) -- (0.67\linewidth,-0.72);
\draw[magenta, densely dashed, thick] (0.75\linewidth,-1.1) -- (0.67\linewidth,-0.72);
\draw[magenta, densely dashed, thick] (0.55\linewidth,-0.72) -- (0.47\linewidth,-0.72);
\draw[magenta, densely dashed, thick] (0.36\linewidth,-0.72) -- (0.28\linewidth,-0.72);
\draw[magenta, densely dashed, thick] (0.13\linewidth,-0.72) -- (0.05\linewidth,-1.25);

	\end{tikzpicture}\vspace{-1.5em}
	\caption{Scheme of combinations per experiment along the four fronts in the E2E approach. For instance, Exp. D labeled with \textcolor{black}{- -} employs the binary $\{-1,1\}$ CA in the SPC optical coding system for the classification task using the regularizer of transmittance level described in \eqref{eq:transmittance_fixed}. }
		\label{fig:combinations}
\end{figure*}
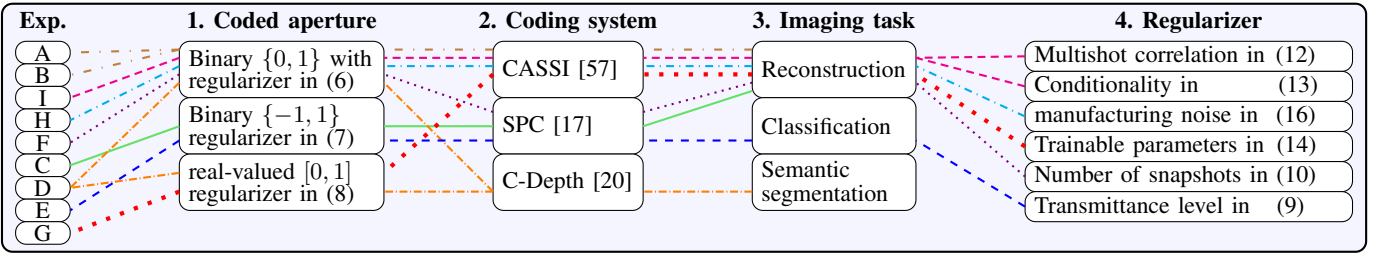

	\subsection{Data Driven Conditionality Constraint}
	\label{sub:conditional}
The conditionality constraint accounts for the sensing matrix structure $\tilde{\mathbf{H}}_{\boldsymbol{\Phi}}$, which should be linearly independent along its rows and columns, to ease its inversion while preserving the features after projection~\cite{hinojosa2018coded}. 
Authors in~\cite{hinojosa2018coded,mejia2018binary,mojica2017high} address this aim
by imposing the regularizer $\| \tilde{\mathbf{H}}_{\boldsymbol{\Phi}}^T\tilde{\mathbf{H}}_{\boldsymbol{\Phi}} -c\mathbf{I} \|_F^2$, where $\mathbf{I}$ denotes the identity matrix for a constant $c\in \mathbb{R_{++}}$. We extend this regularizer by taking advantage of the available data such that the conditionality is improved according to the specific dataset of interest instead of for general acquisition. Hence, our proposal introduces the following regularizer in \eqref{eq:regularization1} to promote a data-driven improved conditionality in the sensing matrix
\begin{equation}
	 R_7(\tilde{\boldsymbol{\Phi}}) = \frac{1}{K} \sum_{k} \left\| \tilde{\mathbf{H}}_{\boldsymbol{\Phi}}^T\tilde{\mathbf{H}}_{\boldsymbol{\Phi}}\mathbf{f}_k -\mathbf{f}_k \right\|_2^2.
	\label{eq:conditionality}
	\end{equation}

	\subsection{Modeling Considerations}
	\label{sub:modeling}
	\subsubsection{Trainable Parameters}
	The number of trainable parameters is a key aspect in the efficiency and performance of the proposed E2E approach. Specifically, when the optical layer contains a set of $S$ 2D or 3D CAs with $N\times M$ and $N\times M \times L$ elements, respectively, the number of trainable parameters ascends to $SMN$ and $SMNL$, respectively. This limits its use for large scale scenes because of the expensive computational memory requirements and the potential over-fitting problems~\cite{goodfellow2016deep}.
We propose to reduce the number of trainable parameters by adding spatial structure to each CA so that a kernel $\mathbf{Q}^s$ of size $\Delta_n \times \Delta_m$ $\ll N \times M$ is periodically repeated as follows
	\begin{equation}
	\boldsymbol{\Phi}^s = \mathbf{1} \otimes \mathbf{Q}^s,
	\label{eq:Kernel_trainable}
	\end{equation} 
	for $s=1,\ldots,S$, where $\mathbf{1}$ denotes a matrix of size $\frac{N}{\Delta_{n}} \times \frac{M}{\Delta_m} $ with all elements equal to $1$, and $\otimes$ represents the Kronecker product. 
	Including such periodicity reduces the total number of trainable variables to $S\Delta_n\Delta_m$ and $S\Delta_n\Delta_mL$ for the 2D and 3D representation, respectively. Note that this model formulation enables to train some optical coding systems with small portions of the training images, known as patches, by training directly $\mathbf{Q}^s$~\cite{gelvez2020nonlocal}. In some 3D applications, this training parameter can be further reduced. For instance, in the colored CA each color pixel can be expressed as a linear combination of $V < L$ fixed optical filters $\{\mathbf{w}^v \in [0,1]^{L}\}_{v=1}^{V}$, so that, each element of the 3D kernel can be rewritten as
	\begin{equation}
	\mathbf{Q}_{i,j,\ell} = \sum_{n} \mathbf{w}_{\ell}^{n}\mathbf{A}_{i,j}^{n},
	\label{eq:Colored_CA}
	\end{equation}
	where $\mathbf{A}$ denotes the trainable weights of the linear combinations, and in consequence, the number of trainable parameters is reduced to $S\Delta_n\Delta_mV$.

		\subsubsection{Manufacturing Noise}
The manufacturing noise represents a problem in the design of optical elements since it can lower the obtained benefits with the design when implemented in a real setup~\cite{diaz2019adaptive,sitzmann2018end}. To overcome this problem, an exhaustive calibration process over the real setup is carried out to achieve a performance as reliable as the sensing model used in the simulations. On the other hand, this problem can also be addressed by refining the noise modeling in the design. Therefore, we aim at considering two sources of perturbations at each step of the forward propagation in the E2E optimization, one into the projected measurements as expressed in \eqref{eq:all_measurments}, which commonly comes from the level of illumination~\cite{he2015hyperspectral}, and the other into the CA, which commonly comes from the manufactured process~\cite{sitzmann2018end}. We add the latter manufacturing noise as follows 
	\begin{equation}
	\tilde{\boldsymbol{\Phi}} = \boldsymbol{\Phi} + \boldsymbol{\eta},
	\label{eq:ManNois}
	\end{equation}
	where $||\boldsymbol{\eta}||_{\infty} << ||\boldsymbol{\Phi}||_{\infty}$, and the distribution of the noise $\boldsymbol{\eta}$ varies according to fabrication processes.
	
Notice that these regularizes have been designed for specific purposes. However, some of them can be incorporated into the optimization problem to get the desired behavior.
\section{Simulation and Results}
		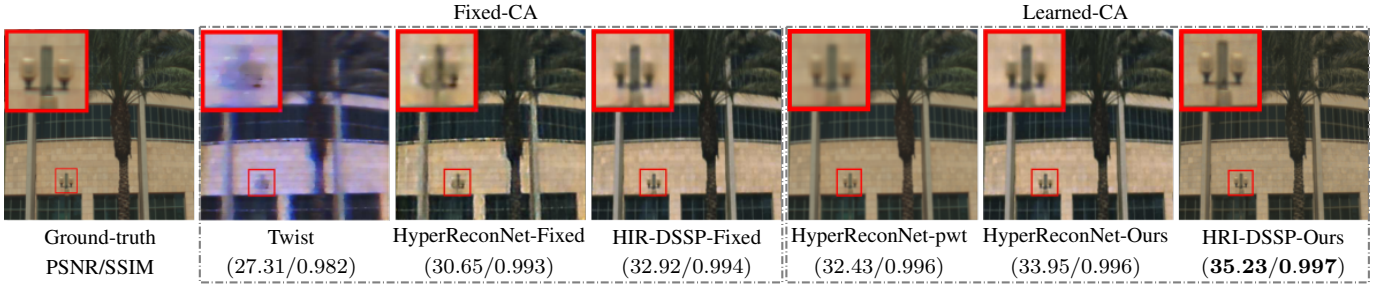
\begin{figure*}[!t]
		\begin{tikzpicture}
		 \tikzstyle{every node}=[font=\footnotesize]
		\node at (0,0){
		\includegraphics[width=1\linewidth]{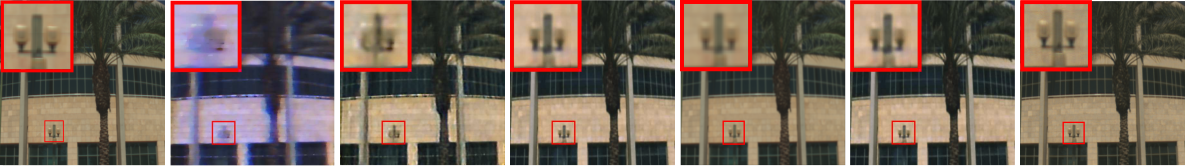}};

		\node at (-0.43\linewidth,-1.5cm) {Ground-truth};
		\node at (-0.43\linewidth,-1.9cm) {PSNR/SSIM};
		
		\draw[gray, densely dashdotted, thick] (-0.357\linewidth, 1.3) rectangle (0.07\linewidth,-2.1);
		\node at (-0.14\linewidth, 1.5cm) {Fixed-CA};
		\node at (-0.29\linewidth,-1.5cm) {Twist};
		\node at (-0.29\linewidth,-1.9cm) {$(27.31/0.982)$};
		\node at (-0.145\linewidth,-1.5cm) {HyperReconNet-Fixed};
		\node at (-0.145\linewidth,-1.9cm) {$(30.65/0.993)$};
		\node at (0.00\linewidth,-1.5cm) {HIR-DSSP-Fixed};
		\node at (0.00\linewidth,-1.9cm) {$(32.92/0.994)$};
		\draw[gray, densely dashdotted, thick] (0.073\linewidth, 1.3) rectangle (0.5\linewidth,-2.1);
		\node at (0.285\linewidth, 1.5cm) {Learned-CA};
		\node at (0.141\linewidth,-1.5cm) {HyperReconNet-pwt};
		\node at (0.141\linewidth,-1.9cm) {$(32.43/0.996)$};
		\node at (0.285\linewidth,-1.5cm) {HyperReconNet-Ours};
		\node at (0.285\linewidth,-1.9cm) {$(33.95/0.996)$};
		\node at (0.43\linewidth,-1.5cm) {HRI-DSSP-Ours};
		\node at (0.43\linewidth,-1.9cm) {$\mathbf{(35.23/0.997)}$};
		\end{tikzpicture}\vspace{-2em}
	 \caption{
	 \textcolor{black}{RGB mapping comparison of the reviewed data-driven approaches, employing fixed and learned CA into the network. Notice that, our CA design regularization improves the quality of previous frameworks for the recovery task.}\vspace{-0.5em}}
	\label{fig:simulation_experiments}
	\end{figure*}
 This section describes the experiments that quantify the performance and validate the effectiveness of the proposed customizable E2E approach for CI tasks. Each experiment considers a combination between the elements of the four fronts mentioned below, and summarized in Fig. \ref{fig:combinations} to ease the reading and understanding. 

 \subsubsection{Three categories of CAs} BCA $\boldsymbol{\Phi}_{i,j}\in \{0,1\}$; BCA $\boldsymbol{\Phi}_{i,j}\in \{-1,1\}$; and real-valued CA with five quantization levels, i.e. $\boldsymbol{\Phi}_{i,j}\in [0,0.25,0.5,0.75,1]$.

 \subsubsection{Three optical coding systems}\footnote{The code with some interactive examples can be found \url{https://github.com/jorgebaccauis/Deep_Coded_Aperture}} The single pixel camera (SPC)~\cite{duarte2008single} using the MNIST dataset~\cite{lecun2010mnist} with 60,000 images of $28 \times 28$ spatial pixels of hand written numbers from 0 to 9, \textcolor{black}{where the input corresponds to different shots according to the compression ratio, and each number denotes a class. A detailed description of the employed SPC testbed can be found in the Supplementary Material.} The coded aperture snapshot spectral imager (CASSI)~\cite{wagadarikar2008single} \textcolor{black}{using two public datasets: the ARAD hyperspectral dataset~\cite{arad2020ntire} with $480$ spectral images of $482 \times 512$ spatial pixels and $31$ spectral bands from $400$nm to $700$nm with a $10$nm step, and the ICVL dataset~\cite{arad2016sparse}, with $201$ images of $1392\times 1300$ spatial pixels and the same spectral resolution of ARAD. For ICVL we obtained 9 overlapping patches of $512 \times 512$ of each image and followed the principles in ~\cite{wang2018hyperreconnet,wang2019hyperspectral} to partition the training and testing sets}. 
 The coded aperture with conventional camera for depth estimation (C-depth)~\cite{levin2007image} \textcolor{black}{was simulated using a $50mm$ $f/1.8$ lens similar to~\cite{levin2007image}, and} using the NYU Depth Dataset with $1449$ RGB images \textcolor{black}{of continuous depth images with depth from 0 to 10 meters captured by Microsoft Kinect~\cite{Silberman:ECCV12}. We use the segmentation label for $40$ classes provided in~\cite{gupta2013perceptual} and the standard  training/test split with $795$ and $654$ images, respectively~\cite{Silberman:ECCV12}. We discretize the depth into $15$ depth maps following the closest value of these depths values  $-3\log(x)+8$  for $x$ from $0.9$ to $11$.}

 \subsubsection{Three CI tasks} reconstruction, classification, and semantic segmentation with different neural networks as explained in each experiment, where we start with the shared weights of the used network, but no layer is frozen. The mean-square error (MSE) and the cross-entropy metrics were used as the loss function for the reconstruction and classification tasks, respectively. On the other hand, the cross-entropy with the softmax activation functions were used as the loss function for the semantic segmentation task. Notice that the E2E approach can be adapted to any pre-existing DNN~\footnote{More details about the networks and their tuning process can be found in Supplementary Material.}.
 
 \subsubsection{Proposed regularizers} functions described in \ref{sec:Design}.

 The metrics to measure the performance of the proposal vary according to the task as follows, the MSE, spectral angle mapper (SAM) and the peak-signal-to-noise ratio (PSNR), calculated as in~\cite{bacca2019noniterative}, are used for reconstruction; the average classification accuracy, calculated as in~\cite{lecun1998gradient}, is used for classification and semantic segmentation; and the mean score from the intersection over union (IoU) calculated as in~\cite{long2015fully}, is used as additional metric for semantic segmentation task.

\subsection{\textcolor{black}{Validation of the proposed regularizers}}
\label{sec:simulated_setup}
\textcolor{black}{
This experiment aims to show that coupling the design of the sensing protocol and the decoder training increase the quality of CI tasks. For this, we compared the results against the non data-driven recovery method TwIST with TV prior~\cite{kittle2010multiframe}, and the data-driven recovery networks Hyperspectral Image Reconstruction using a Deep Spatial-Spectral Prior (HIR-DSSP)~\cite{wang2019hyperspectral}, and HyperReconNet~\cite{wang2018hyperreconnet}. For Twist, we employed a fixed CA. For (HIR-DSSP) and HyperReconNet, we evaluated different variations: first, we evaluated the performance when using a fixed CA generated following a Bernoulli distribution with parameter $0.5$, whose results are denoted as HyperReconNet-Fixed, and HIR-DSSP-Fixed. Furthermore, we evaluated the performance when joining the design of the CA and the network weights through our proposed binary regularization, whose results are denoted as HyperReconNet-Ours, and HIR-DSSP-Ours. Finally, we evaluated the performance when using the strategy proposed by the HyperReconNet~\cite{wang2018hyperreconnet}, which learns the CA using a piece-wise threshold (pwt), whose result is denoted as HyperReconNet-pwt. Figure \ref{fig:simulation_experiments} shows an RGB visual representation of the recovered images with the PSNR and SSIM quantitative metrics which show that the proposed method outperforms state-of-the-art methods. \textcolor{black}{We remark} the versatility of the proposed regularizers which can be straightforward incorporated at any pre-existed net.}

\subsection{Validation in a real setup experiment}

\begin{figure}[b!]
		\centering
		\includegraphics[width=0.98\linewidth]{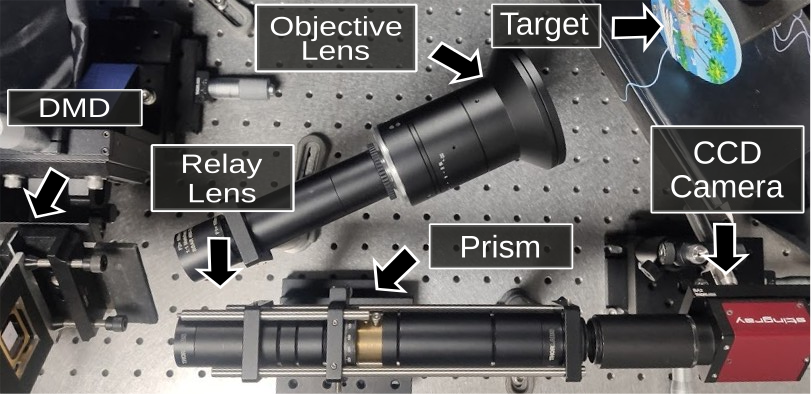}\vspace{-0.6em}
		\caption{Testbed CASSI implementation where the relay lens focuses the encoded light by the DMD into the sensor after dispersed by the prism.}	\label{fig:realsetup}
	\end{figure}

\textcolor{black}{Section \ref{sec:simulated_setup}} demonstrated that coupling the sensing protocol design and the processing method increases the quality of CI tasks; however, most of the obtained benefits are lowered when applying those methods in real setups~\cite{correa2016multiple,bacca2019noniterative}. Hence, this experiment validates the E2E approach with a real setup corresponding to one single snapshot of the CASSI testbed laboratory implementation depicted in Fig. \ref{fig:realsetup} 
\textcolor{black}{For this, the ARAD dataset was used to train the proposed approach with a compression ratio $\gamma= 0.034$.} The setup contains a $100$-$nm$ objective lens, a high-speed digital micro-mirror device (DMD) (Texas Instruments-DLI$4130$), with a pixel size of $13.6 \mu m$, an Amici Prism (Shanghai Optics), and a CCD (AVT Stingray F-$145B$) camera with spatial resolution $1388 \times 1038$, and pitch size of $6.8 \mu m$. \textcolor{black}{The $482\times 512$ designed CA is placed at center of the DMD.}

 The performance was compared against the results obtained with the same DNN and hyper-parameter tuning process, \textcolor{black}{but with a fixed CA and learning only the network weights,} i.e., the main difference is that the \textcolor{black}{optical layer} is trainable in the proposed method. For the fixed CA, we used a random CA and a designed blue-noise CA~\cite{correa2016spatiotemporal}. 
 \begin{figure}[t!]
		\begin{tikzpicture}
		 \tikzstyle{every node}=[font=\small]
		\node at (0,0.33) {
		\includegraphics[width=0.48\textwidth]{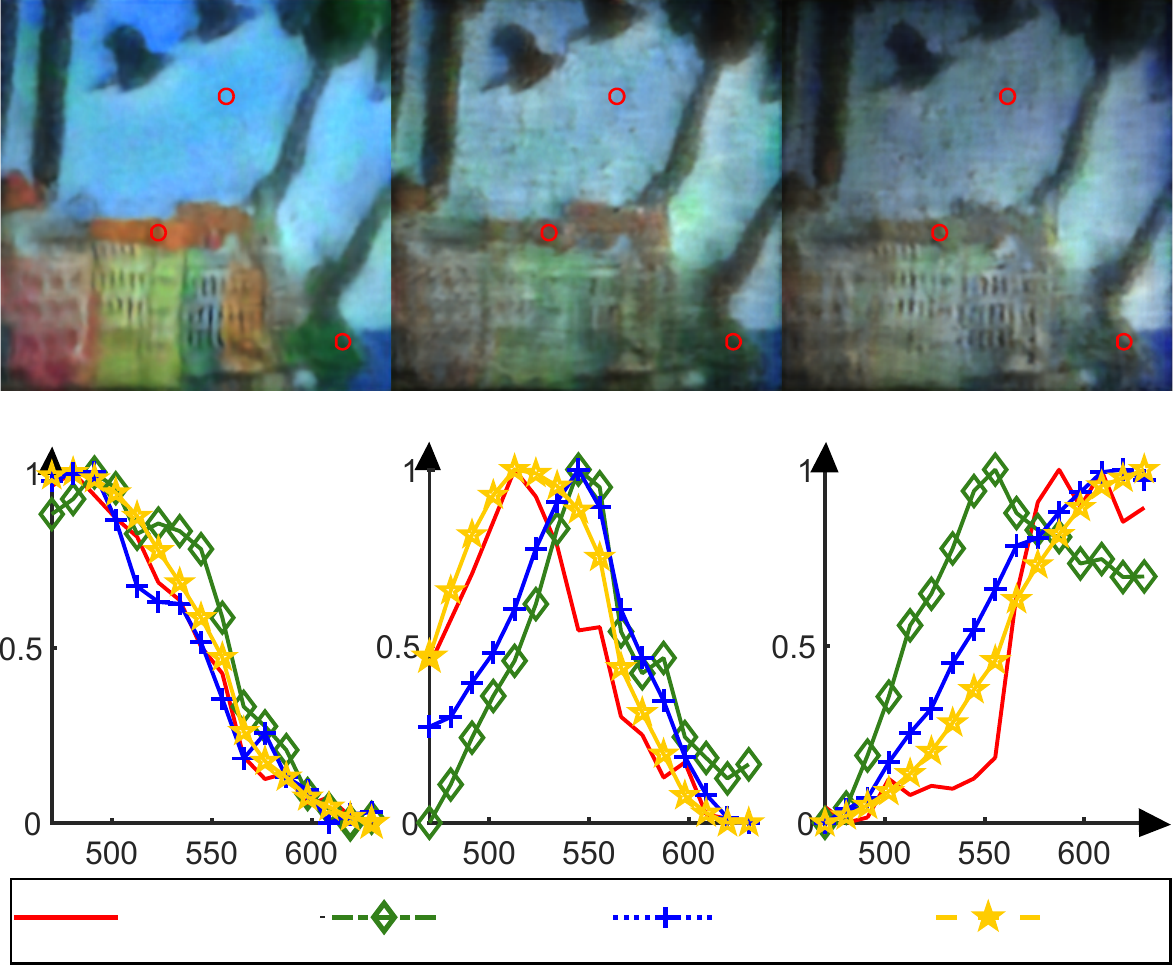}};
		\node at (-0.3333\linewidth,4.1cm) {DNN-Design};
		\node at (0\linewidth,4.1cm) {Blue-noise};
		\node at (0.3333\linewidth,4.1cm) {Random};
		\draw[densely dashdotted,thick] (-0.5\linewidth,0.92) -- (0.5\linewidth,0.92);
		\node at (-0.3333\linewidth,0.71cm) {$\mathbf{P_1}$};
		\node at (0\linewidth,0.71cm) {$\mathbf{P_2}$};
	 \node at (0.3333\linewidth,0.71cm) {$\mathbf{P_3}$};
	 \node at (-0.17\linewidth,-2.4cm) {$\lambda$};
	 \node at (0.14\linewidth,-2.4cm) {$\lambda$};
	 \node at (0.46\linewidth,-2.4cm) {$\lambda$};
	 \node[anchor=west] at (-0.4\linewidth,-2.9cm) {Reference};
	 \node[anchor=west] at (-0.135\linewidth,-2.9cm) {Random};
	 \node[anchor=west] at (0.095\linewidth,-2.9cm) {Blue-noise};
	 \node[anchor=west] at (0.37\linewidth,-2.9cm) {DNN};
	 
\node[red] at (-0.33\linewidth,3.5cm) {$\mathbf{P_1}$};
\node[red] at (-0.39\linewidth,2.5cm) {$\mathbf{P_3}$};
\node[red] at (-0.24\linewidth,1.7cm) {$\mathbf{P_2}$};
\node[red] at (-0.01\linewidth,3.5cm) {$\mathbf{P_1}$};
\node[red] at (-0.07\linewidth,2.5cm) {$\mathbf{P_3}$};
\node[red] at (0.08\linewidth,1.7cm) {$\mathbf{P_2}$};
\node[red] at (0.31\linewidth,3.5cm) {$\mathbf{P_1}$};
\node[red] at (0.25\linewidth,2.5cm) {$\mathbf{P_3}$};
\node[red] at (0.4\linewidth,1.7cm) {$\mathbf{P_2}$};

\node at (-0.23\linewidth,0.43cm){\footnotesize{SAM}};
\draw (-0.275\linewidth,0.25cm) rectangle (-0.185\linewidth,-0.65cm);
\node[my_green] at (-0.23\linewidth,0.1cm) {\footnotesize{$0.172$}};
\node[blue] at (-0.23\linewidth,-0.2cm) {\footnotesize{$0.102$}};
\node[my_yellow] at (-0.23\linewidth,-0.5cm) {\footnotesize{$0.054$}};

\node at (0.1\linewidth,0.43cm){\footnotesize{SAM}};
\draw (0.055\linewidth,0.24cm) rectangle (0.145\linewidth,-0.65cm);
\node[my_green] at (0.1\linewidth,0.1cm) {\footnotesize{$0.648$}};
\node[blue] at (0.1\linewidth,-0.2cm) {\footnotesize{$0.463$}};
\node[my_yellow] at (0.1\linewidth,-0.5cm) {\footnotesize{$0.146$}};

\node at (0.41\linewidth,-0.65cm){\footnotesize{SAM}};
\draw (0.369\linewidth,-0.85cm) rectangle (0.455\linewidth,-1.75cm);
\node[my_green] at (0.41\linewidth,-1.0cm) {\footnotesize{$0.630$}};
\node[blue] at (0.41\linewidth,-1.3cm) {\footnotesize{$0.316$}};
\node[my_yellow] at (0.41\linewidth,-1.6cm) {\footnotesize{$0.224$}};
	\end{tikzpicture}\vspace{-1.7em}
		\caption{(Top) RGB visual representation of the three evaluated methods (Net-design, Blue-noise and random). (Bottom) Comparison of the normalized spectral signatures at three points in the recovered scenes. At each point it is shown the quantitative SAM metric to quantify the improvement.}
		\label{fig:real_measurments}\vspace{-0.5em}
	\end{figure}
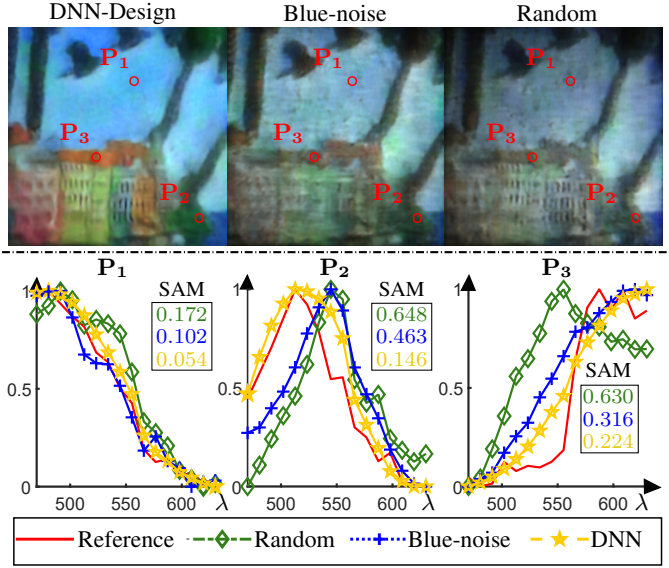

 Figure \ref{fig:real_measurments} (Top) illustrates an RGB visual representation comparison of the spatial quality between the reconstructions obtained along the three approaches. It can be noticed a significant improvement in the visual results for a real setup when coupling the design of the CA in the DNN. Fig. \ref{fig:real_measurments} (Bottom) illustrates a comparison of the spectral quality at three random spatial locations whose spectral response was measured in the laboratory with the commercially available spectrometer (Ocean Optics USB2000+). It can be noticed that the coupled approach decreases the spectral angle between the estimated and reference spectral signatures in comparison to the fixed CAs. 



\subsection{Implementation constraints experiment}
    The regularization strategy to customize the CA properties formulated in \eqref{eq:regularization1} involves the selection of a proper regularization parameter $\rho_q$. This experiment aims to show the importance of such selection to obtain targeted quantization levels in the CA using regularizer $R_2(\cdot)$ in \eqref{eq:family_gray_0}, which indirectly leads to an easier implementation of the CA, and how it implies a trade-off with the loss function, responsible for the task quality.

	
	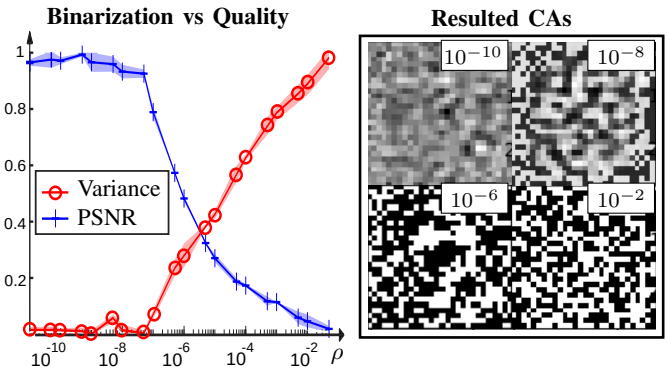
\begin{figure}[b!]
		\begin{tikzpicture}
		 		 \tikzstyle{every node}=[font=\small]
		 		\node[inner sep=0pt] (Static) at (0\linewidth,0){
	 \includegraphics[width=0.99\linewidth]{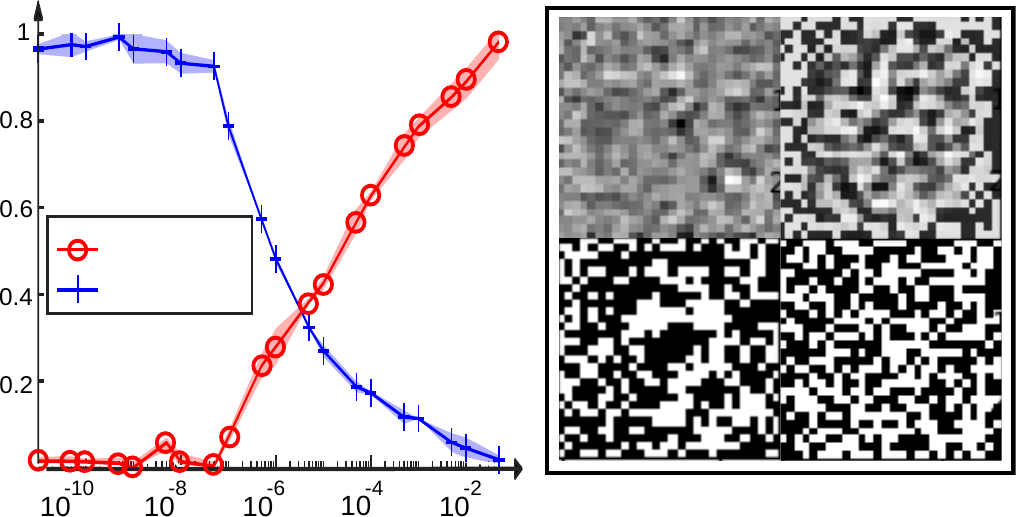}};
	 \node at (-0.25\linewidth,2.4) {\textbf{Binarization vs Quality}};
	 \node at (0.25\linewidth,2.4) {\textbf{Resulted CAs}};
	 \node[anchor=west] at (-0.4\linewidth,0.14) {Variance};
	 \node[anchor=west] at (-0.4\linewidth,-0.25) {PSNR};
	 \node[anchor=west] at (-0.02\linewidth,-2.1) {$\rho$};
	 
	 \filldraw [fill=white, draw=black] (0.16\linewidth,1.7) rectangle (0.26\linewidth,2.1);
 \node[anchor=west] at (0.15\linewidth,1.9) {\footnotesize{$10^{-10}$}};
 
 	 \filldraw [fill=white, draw=black] (0.38\linewidth,1.7) rectangle (0.48\linewidth,2.1);
 \node[anchor=west] at (0.38\linewidth,1.9) {\footnotesize{${10^{-8}}$}};

 	 \filldraw [fill=white, draw=black] (0.16\linewidth,-0.2) rectangle (0.26\linewidth,0.2);
 \node[anchor=west] at (0.16\linewidth,-0.0) {\footnotesize{$10^{-6}$}};
 
 	 \filldraw [fill=white, draw=black] (0.38\linewidth,-0.2) rectangle (0.48\linewidth,0.2);
 \node[anchor=west] at (0.38\linewidth,-0.0) {\footnotesize{$10^{-2}$}};
 
	 \end{tikzpicture}\vspace{-1.0em}
	 \caption{trade-off between the implementability in terms of normalized variance in logarithmic scale, and the obtained quality in terms of normalized PNSR.}
	 \label{fig:QI-trade-off}
	\end{figure}

	 Figure \ref{fig:QI-trade-off} illustrates this trade-off where the achievement of the discretization levels is measured in terms of the logarithmic variance of the elements in the CA,
	 and the quality is measured in terms of the PNSR, along an equally spaced log scale interval of $\rho$, from $1e^{-11}$ to $1e^{-2}$. For comparison and analysis purposes the resulting PSNR and variances were normalized along the experiments. It can be observed that the highest the value of $\rho$, the highest the achieved binarization implying the easiest implementable CA at the cost of the poorest task quality. Unlike, the smallest the value of $\rho$, the highest the task quality at the cost of non-full binarization indicating the hardest implementable CA.
	 

	To mitigate this trade-off effect, we present a dynamic strategy for $\rho$ that aims to start with a very small value, $\rho^0$, and increase its value a factor $\alpha >1$ every $\beta \leq T$ epochs, in such a manner that a target value $\rho^T$ is achieved after executing $T$ epochs. This strategy allows the network to fit and benefit the particular task during the first epochs when $\rho$ is small, and then guarantee the binarization of the CA during the last epochs when $\rho$ is large. In this manner, $\rho$ is updated a total number of $\lfloor\frac{T}{\beta}\rfloor$ times by using the following update equation
	\begin{equation}
	 \rho^{k+1} = \alpha \rho^{k} ,
	\end{equation}
	for $k=1, \ldots, \lfloor\frac{T}{\beta}\rfloor$.
	Further, notice that this multiplicative update strategy leads to the relation $ 
	 \rho^T = \alpha^{\lfloor\frac{T}{\beta}\rfloor}\rho^0,$
	where, $\rho^0<\rho^T \in \mathbb{R}^+, \beta\in [1,T]
	,$ and $T \in [1,\infty]$ can be seen as hyperparameters and $\alpha$ can be calculated as
	
	\begin{equation}
	 \alpha = \sqrt[{\lfloor\frac{T}{\beta}\rfloor}]{{\rho^{T}/\rho^{0}}}. 
	 \label{eq:alpha}
	\end{equation}
	
	\begin{figure}[t!]
 \begin{tikzpicture}
		 		 \tikzstyle{every node}=[font=\small]
		 		\node[inner sep=0pt] (Static) at (0\linewidth,0){
 \includegraphics[width=1\linewidth,left]{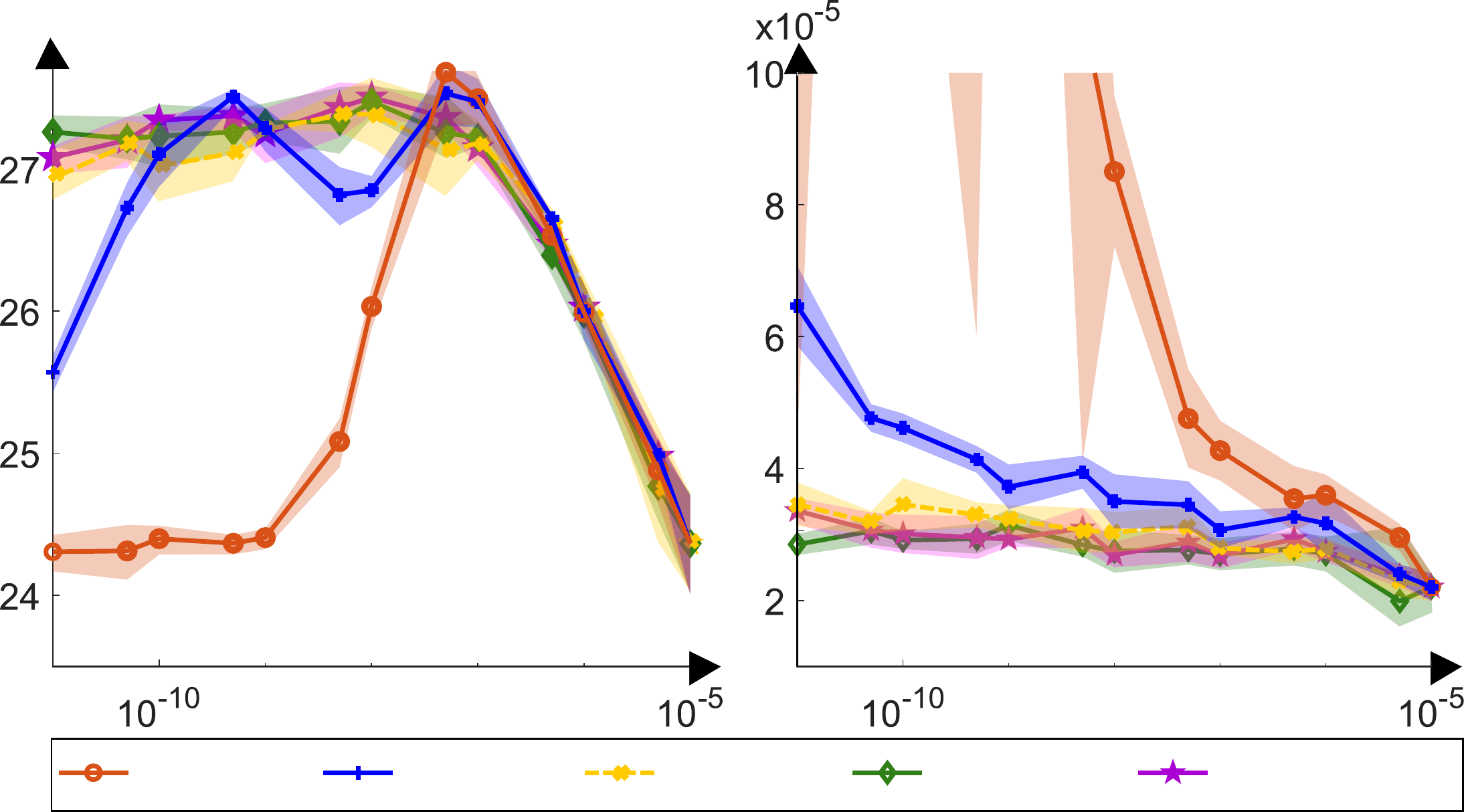}}; 
 
\node at (-0.345\linewidth,-2.25) {$\lfloor\frac{T}{\beta}\rfloor = 1$};
\node at (-0.17\linewidth,-2.25) {$\lfloor\frac{T}{\beta}\rfloor = 2$};
\node at (0.01\linewidth,-2.25) {$\lfloor\frac{T}{\beta}\rfloor = 5$};
\node at (0.2\linewidth,-2.25) {$\lfloor\frac{T}{\beta}\rfloor = 20$};
\node at (0.41\linewidth,-2.25) {$\lfloor\frac{T}{\beta}\rfloor = 100$};
 	 \node[anchor=west] at (0.175\linewidth,2.4) {\textbf{Variance}};
 	 \node[anchor=west] at (0.46\linewidth,-1.3) {$\rho^0$};
 	 \node[anchor=west] at (-0.04\linewidth,-1.3) {$\rho^0$};
 	 \node[anchor=west] at (-0.33\linewidth,2.4) {\textbf{Quality in PSNR}};
 	 \node at (-0.46\linewidth,2.4) {[db]};
 \end{tikzpicture}\vspace{-1.6em}
 \caption{Depiction of the dynamic $\rho$ strategy along five cases of the ratio $\lfloor\frac{T}{\beta}\rfloor$. (Left) quality measured in terms of PSNR [db], (right) binarization level in the CA measured in terms of the variance of its entries.}
 \label{fig:DynamicRho}
\end{figure}
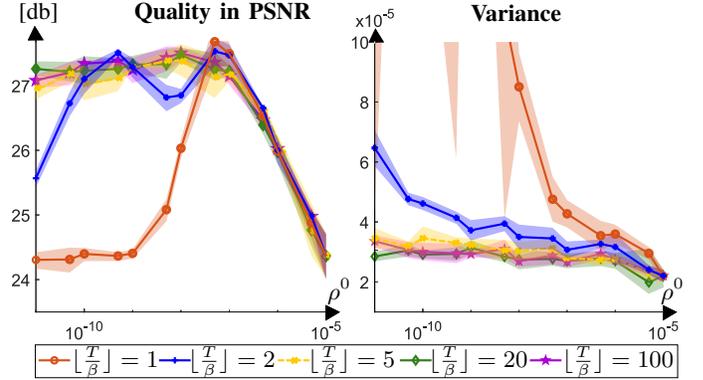
	
	Figure \ref{fig:DynamicRho} shows the behaviour of the dynamic $\rho$ strategy, for $T=100$ and $\rho^T=1e^{-5}$, along an equally spaced log scale interval of the initial value $\rho^0$ from $1e^{-5}$ to $1e^{-11}$, when updating the parameter $\lfloor\frac{T}{\beta}\rfloor = [1 \hspace{2mm} 5 \hspace{2mm} 20 \hspace{2mm} 100]$ times. The value of $\alpha$ is calculated for each experiment using \eqref{eq:alpha}. It can be seen an improvement in up to 3[db] when using the dynamic strategy in comparison to the quality of the baseline of a static parameter $\rho=1e^{-5}$, while achieving the same binarization determined by the variance of the CA entries. Further, it is observed that when the parameter $\rho$ is updated only once, i.e $\lfloor\frac{T}{\beta}\rfloor = 1$, the behaviour of the dynamic $\rho$ is not stable. \textcolor{black}{We select the value of $\beta$ such that depending on the initial $\rho^0$ and $\alpha=10$, is updated $10$ times along the established number of epochs $T$ in the following of the experiments.}

\subsection{real-valued CA and BCA comparison experiment}

	This experiment compares the performance of the real-valued CA against the BCA for the semantic segmentation task using a C-depth system. To the best of our knowledge, this is the first work that employs the C-Depth optical system for this high-level task. \textcolor{black}{Hence, for comparison purposes we evaluate the performance for the cases when employing only RGB and RGB+Depth (RGB+D) information.
		 Specifically, we evaluate the RefineNet~\cite{lin2017refinenet} that only employs the RGB image as input, the RDFNet~\cite{park2017rdfnet} that extends the RefineNet architecture for RGB+D, and REDFNet that corresponds to the proposed modification of the RDFNet,  such that, it receives an RGB with encoded depth as blurred versions at each depth map,  as the input (See Supplementary Material for more information). In addition, we incorporated comparisons against the state-of-the-art Geometry-Aware~\cite{jiao2019geometry} that employs RGB + estimated D, understanding RGB + estimated D as their strategy of employing depth information only during the training process, such that, during the inference process it only requires the RGB image}. The quantitative quality of the semantic segmentation results measured in terms of accuracy, and IoU are summarized in Table~\ref{tab:sematic_depth}.
		 \begin{table}[t!]
	\scriptsize
	\renewcommand{\arraystretch}{1.1}
	\caption{results using NYUDv2
40-class}
	\label{tab:sematic_depth}
	\centering
	\resizebox{1\columnwidth}{!}{%
	\setlength\tabcolsep{0.05cm}
	\begin{tabular}{llcccc} 
		\toprule 
		Method & Data &Pixel Acc & Mean Acc & Mean IoU \\
		\midrule 
		REDFNet &C-depth BCA & 74.8 & 60.1 & 47.3 \\ 
	 	REDFNet& C-depth real-valued CA &76.0 & 61.4 & 49.9 \\ 
		\midrule 
		RefineNet & RGB & 73.6 & 58.9 & 46.5 \\ 
		RDFNet & RGB+D & 76.0 & 62.8 & 50.1 \\ 
		GA & RGB+estimated D & 84.8 & 68.7 & 59.6 \\ 
		\bottomrule 
	\end{tabular} \vspace{-0.7em} }
\end{table}
 \begin{figure}[t!]
 	\begin{tikzpicture}
 	\tikzstyle{every node}=[font=\footnotesize]
 	\node at (0.2,0.0) {
 		\includegraphics[width=0.92\linewidth]{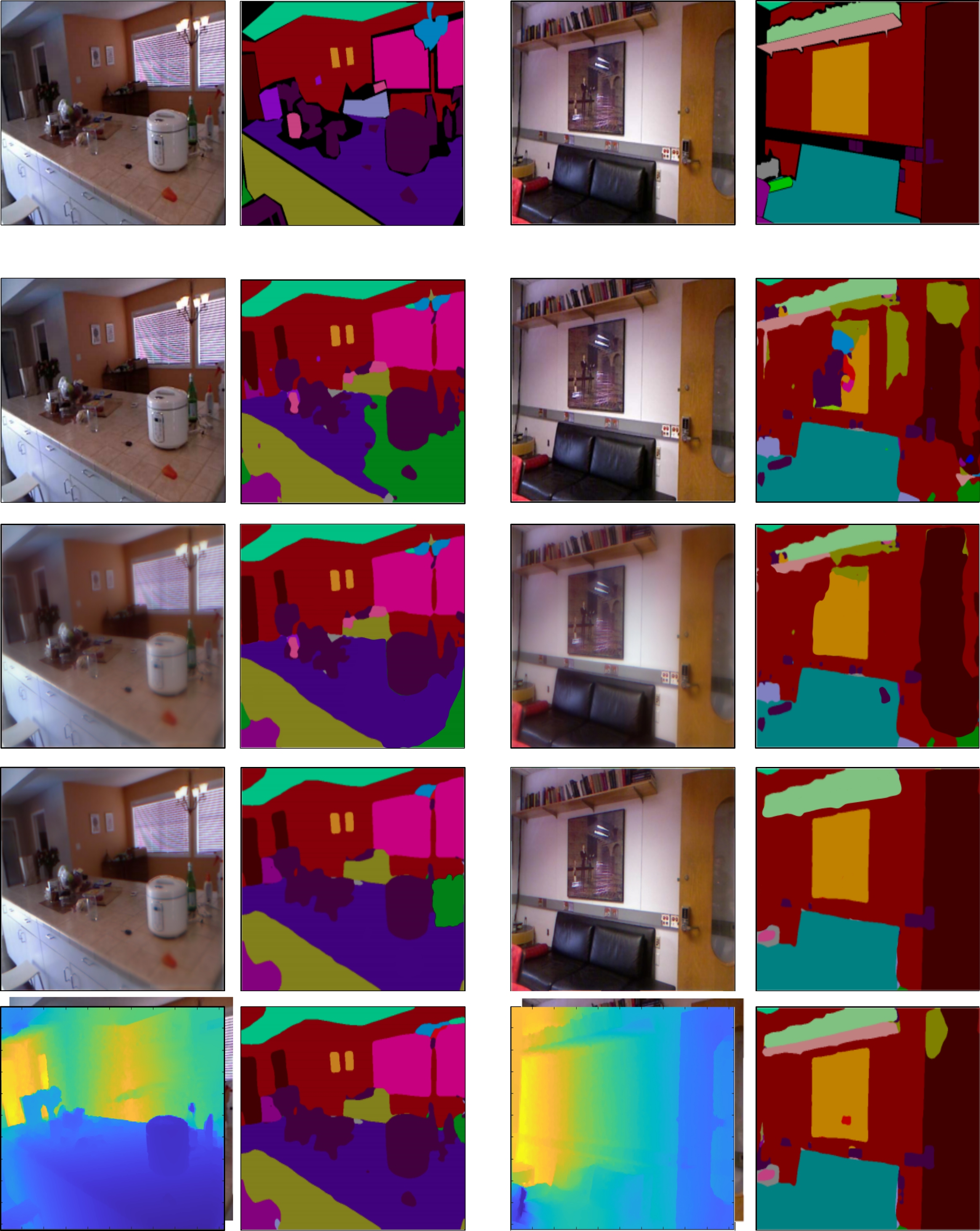}};
 	
 	\node at (-0.22\linewidth,5.7cm) {\textbf{Image 1}};
 	\node at (0.25\linewidth,5.7cm) {\textbf{Image 2}};
 	
 	\node at (-0.33\linewidth,5.3cm) {All-focus};
 	\node at (-0.11\linewidth,5.3cm) {Ground-truth};
 	\node at (0.140\linewidth,5.3cm) {All-focus};
 	\node at (0.37\linewidth,5.3cm) {Ground-truth};
 	
 	\draw[densely dashdotted,thick] (-0.49\linewidth,3.20) -- (0.49\linewidth,3.20);
 	\node at (-0.330\linewidth,3.00cm) {Measurements};
 	\node at (-0.11\linewidth,3.00cm) {Results};
 	
 	\node at (0.14\linewidth,3.00cm) {Measurements};
 	\node at (0.365\linewidth,3.00cm) {Results};

\node[rotate=90] at (-0.5\linewidth,1.9cm) {RefineNet};
\node[rotate=90] at (-0.47\linewidth,1.9cm) {RGB};

\node[rotate=90] at (-0.5\linewidth,-0.15cm) {REDFNet};
\node[rotate=90] at (-0.47\linewidth,-0.15cm) {C-depth BCA}; 
 	
\node[rotate=90] at (-0.50\linewidth,-2.2cm) {REDFNet};
\node[rotate=90] at (-0.47\linewidth,-2.2cm) {real-valued CA};
 	
\node[rotate=90] at (-0.5\linewidth,-4.2cm) {RDFNet};
\node[rotate=90] at (-0.47\linewidth,-4.2cm) {RGB+D}; 	
 	
 	\end{tikzpicture}\vspace{-1.7em}
 	\caption{Visual comparison of semantic segmentation results for two testing images. (Top) All-focus and ground-truth images. (Bottom) Measurements and results obtained with RGB, C-depth BCA, real-valued CA, and RGB+D.}
 	\label{fig:sematic_seg}
 \end{figure}
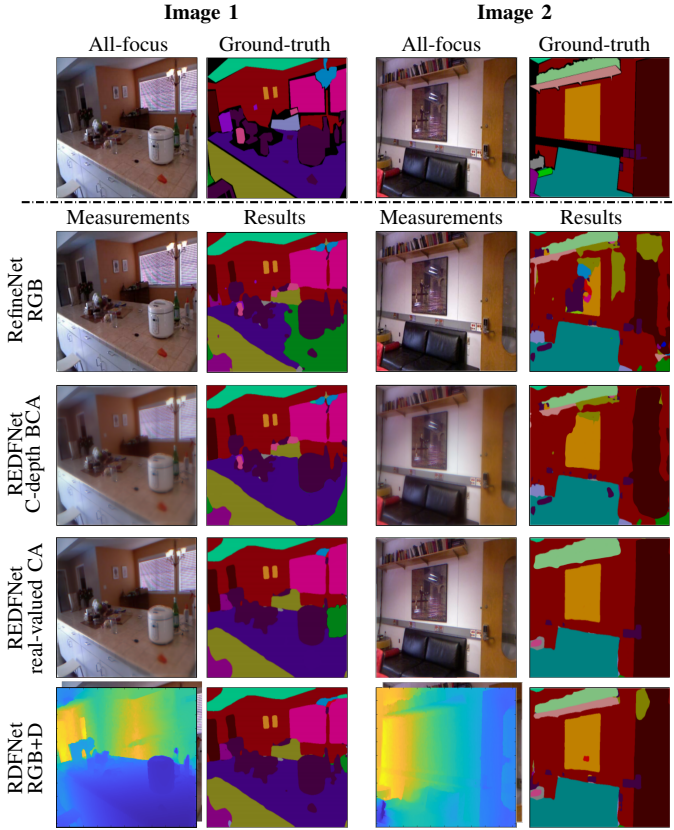
	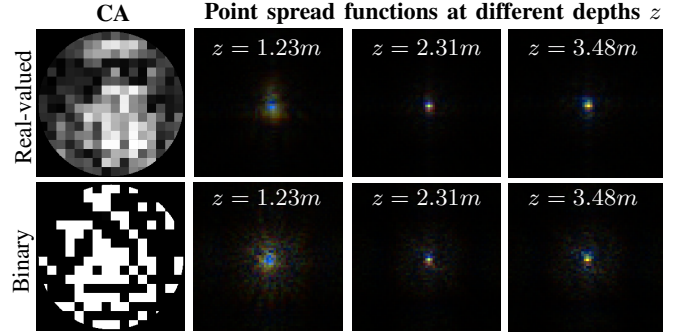
\begin{figure}[t!]
 			\begin{tikzpicture}
 			
 			\tikzstyle{every node}=[font=\small]
 			\node at (0,0.33) {
 		\includegraphics[width=0.94\linewidth]{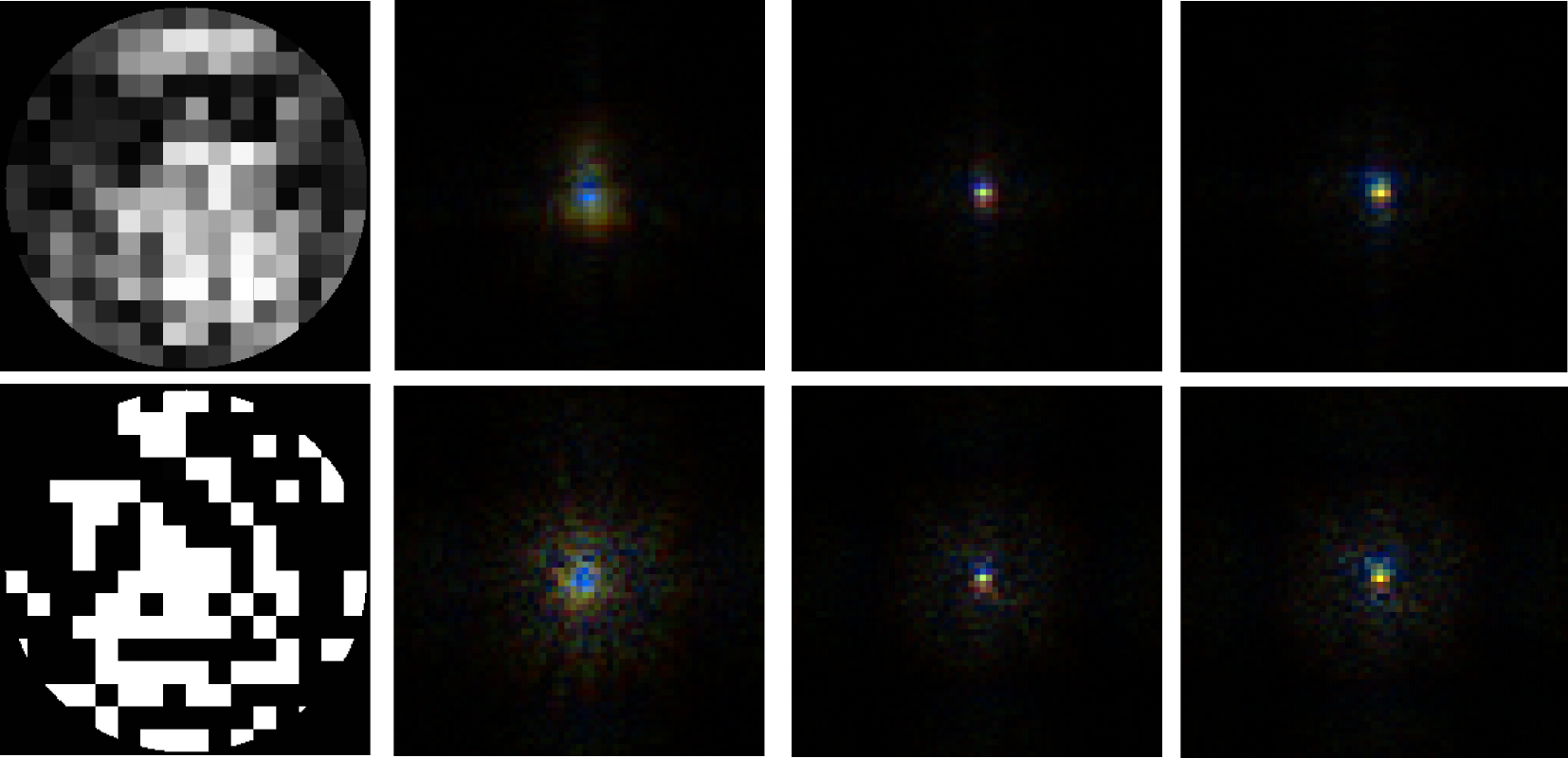}};

 				\node at (-0.355\linewidth,2.55cm) {\textbf{CA}};
 				\node at (0.125\linewidth,2.55cm) {\textbf{Point spread functions at different depths $z$}};

 	\node at (-0.125\linewidth,2.10cm) {\textcolor{white}{$z = 1.23 m$}};
 	\node at (-0.125\linewidth,0.1cm) {\textcolor{white}{$z = 1.23 m$}};
 				
	\node at (0.115\linewidth,2.10cm) {\textcolor{white}{$z = 2.31 m$}};
 \node at (0.115\linewidth,0.1cm) {\textcolor{white}{$z = 2.31 m$}};
 				
	\node at (0.35\linewidth,2.10cm) {\textcolor{white}{$z = 3.48 m$}};
\node at (0.35\linewidth,0.1cm) {\textcolor{white}{$z = 3.48 m$}};

 	\node[rotate=90] at (-0.49\linewidth,1.3cm) {Real-valued};
 	\node[rotate=90] at (-0.49\linewidth,-0.75cm) {Binary};

 			\end{tikzpicture}
 	\centering\vspace{-1.2em}
 	\caption{Visual representation of real-valued CA and BCA, and three point spread function at different $z$, obtained with the proposed E2E approach.}
 	\label{fig:CA_results}
 	
 \end{figure}
	\textcolor{black}{\textcolor{black}{It can be observed that both the real-valued CA and BCA achieve remarkable comparable results concerning the RGB+D setup with the advantage of using only one sensor. Furthermore, the proposed Depth-CA outperforms the semantics results obtained with only an RGB image, which shows the advantage of the depth information for the semantic segmentation task.} 
	\textcolor{black}{These results indicate that the depth information is properly encoded in the projected measurements with the design of the CA into the E2E approach.} \textcolor{black}{Furthermore, the real-valued CA narrowly outperformed the BCA in up to 2\% of pixel accuracy.} Figure \ref{fig:sematic_seg} shows a qualitative comparison of the semantic segmentation results for two testing images. It can be seen the effectiveness of the C-depth system to perform this CI task.}
	 Finally, Fig. \ref{fig:CA_results} depicts the designed real-valued CA and BCA obtained with the proposed E2E approach and their corresponding point spread function (PSF) at different distance, which have a direct relationship with the sensing matrix $\mathbf{H}_{\boldsymbol{\Phi}}$ (See Supplementary Material \textcolor{black}{for a detailed experiment of the obtained CA)}.



\subsection{Transmittance experiment}
The proposed E2E approach provides two options to adjust the transmittance level of the designed CA, (\romannumeral1) by varying the values of $p_d$ in the family of functions in \eqref{eq:family_gray_0}, or (\romannumeral 2) by imposing the transmittance using the regularizer in \eqref{eq:transmittance_fixed}. This experiment evaluates the effect of using these options, separately and together to attain a targeted transmittance along an equally spaced interval from $T_r = 0.1$ to $T_r = 0.9$, using the SPC and for classification task. \textcolor{black}{For that experiment, the initial $\rho$ values of the binary and transmittance regularizer are established as $1e^{-9}$ and $1e^{-7}$, respectively.} To attain each specific transmittance level using the family function option, the values of parameters $p_1$ and $p_2$ were found experimentally as shown in Table II; i.e. to attain $T_r = 0.1$, the found proper values are $p_1=1.8$ and $p_2=1$, and so on. \textcolor{black}{Remark that the family of functions for binary regularization for $ p_1 = 1 $ and $ p_2 = 1 $ converges to the state-of-the-art method in~\cite{bacca2020coupled}.}
\begin{table}[t!]
\begin{center}
\caption{Experimental values for $Tr$, $p_1$, and $p_2$}
	\resizebox{8.7cm}{!} {
\begin{tabular}{|c||c|c|c|c|c|c|c|c|c|}
\hline
\toprule 
\textbf{$T_r$} & 0.1 & 0.2 & 0.3 & 0.4 & 0.5 & 0.6 & 0.7 & 0.8 & 0.9 \\ \hline
\toprule 
\textbf{$p_1$} & 1.8 & 1.5 & 1.3 & 1.125 & 1 & 1 & 1 & 1 & 1 \\ \midrule
\textbf{$p_2$} & 1 & 1 & 1 & 1 & 1 & 1.125 & 1.3 & 1.5 & 1.8 \\ \hline
\toprule 
\end{tabular}}
\label{tab:table1}
\end{center}
\vspace{-1em}
\end{table}
	\begin{figure}[t!]
 	\begin{tikzpicture}
	 \tikzstyle{every node}=[font=\small]
		 		\node[inner sep=0pt] (Examples) at (0\linewidth,0){
		\includegraphics[width=0.96\linewidth,center]{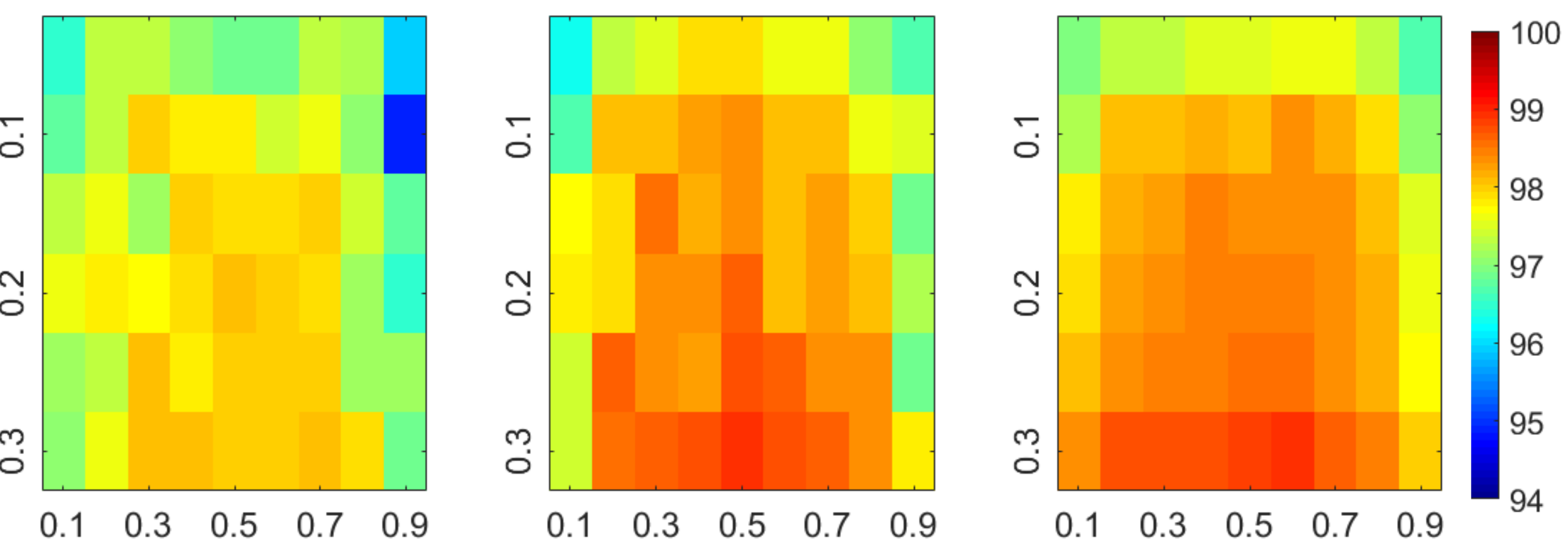}};
		 \node at (0\linewidth,2.0cm) {\textbf{Accuracy map as a function of $Tr$ and $\gamma$ }};
		 \node at (-0.35\linewidth,1.6cm) {{a. Family}};
		 \node at (-0.03\linewidth,1.6cm) {{b. Regularizer}};
		 \node at (0.27\linewidth,1.6cm) {{c. Both}};
		 \node at (-0.3333\linewidth,-1.7cm) {$T_r$};
		 \node at (-0.02\linewidth,-1.7cm) {$T_r$};
		 \node at (0.29\linewidth,-1.7cm) {$T_r$};
		\node[rotate=90] at (-0.5\linewidth,0cm) {$\gamma$};
		\node[rotate=90] at (-0.19\linewidth,0cm) {$\gamma$};
		\node[rotate=90] at (0.12\linewidth,0cm) {$\gamma$};
		\end{tikzpicture}\vspace{-2em}
		\caption{Accuracy map as a function of the transmittance and compression ratio when using a. the family of functions, b. the regularizer, and c. both together.}
		\label{fig:transmitance}
	\end{figure}
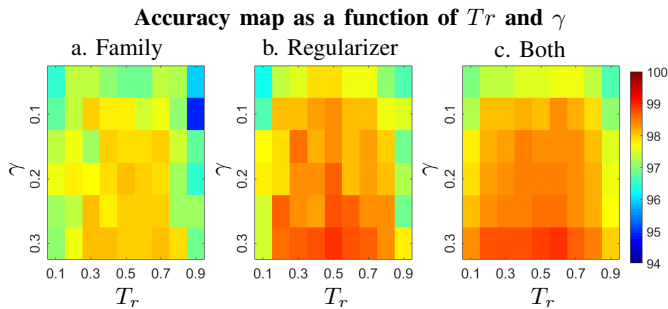
Figure \ref{fig:transmitance} shows a map of the classification accuracy in the scale $[94-100]$ obtained when adjusting the targeted transmittance using a. only the family of functions, b. only the regularizer, and c. both options together, along an equally spaced interval of compression ratios from $0.05$ to $0.3$. It can be seen that the maximum accuracy for each compression ratio is achieved at different transmittance levels. Also, notice that using the two options together provides a better quality for very small and high transmittance levels, which is the desired behavior in the CI applications described in \ref{sub:trans}.

	\subsection{Number of snapshots}
	
		\begin{figure}[b!]
		\begin{tikzpicture}
		 \tikzstyle{every node}=[font=\small]

 \node[inner sep=0pt] (Static) at (0\linewidth,0){
	 \includegraphics[width=1\linewidth]{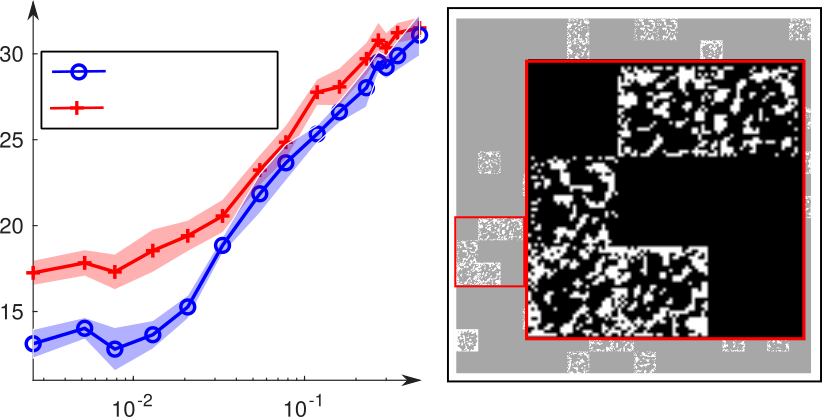}};

	 \node at (-0.22\linewidth,2.4) {\textbf{ Quality in PSNR}};
	 \node at (0.25\linewidth,2.4) {\textbf{Resulting CAs }};
	
	 \node[anchor=west] at (-0.37\linewidth,1.5 ) {Fixed};
	 \node[anchor=west] at (-0.37\linewidth, 1.1) {Proposed};
	 \node[anchor=west] at (-0.04\linewidth,-2.1) {$\gamma$};
	 \node[anchor=west] at (-0.5\linewidth,2.35) {[db]};
	 
	 \end{tikzpicture}\vspace{-1.8em}
	 \caption{Learning the number of snapshots with regularizer in \eqref{eq:regu_shots}. \textit{(left)} Quality in terms of PNSR vs the compression ratio ($\gamma$), for learned and fixed number of snapshots. \textit{(right)} Visual representation of some obtained CAs, where full-zero CAs indicate no acquiring the snapshot.}
	 \label{fig:trade}
	\end{figure}
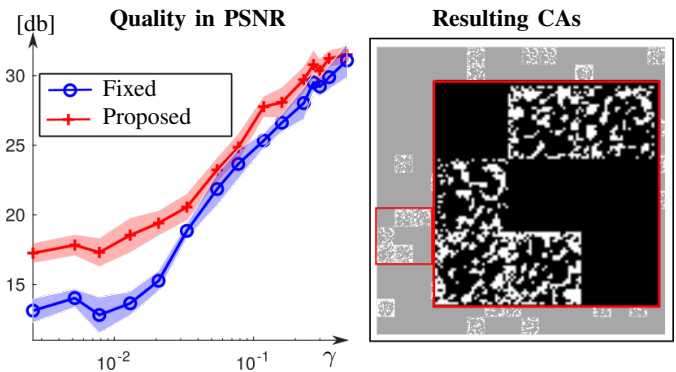

	This experiment evaluates the effectiveness of the proposed regularizer in \eqref{eq:regu_shots} to determine the optimal number of snapshots to achieve an acceptable task performance. 
	For this, we established a maximum number of snapshots of $S' =392$, which is equivalent to a compression ratio of $\gamma = 0.5$. 
	
	\textbf{Numerical Note:} This experiment uses both the binarization and number of snapshots regularizers, which fundamentally pursue opposite purposes, and in consequence, it may result in the number of snapshots not being reduced. Therefore, to guarantee obtaining some CAs with complete zero entries, we modify the forward propagation of $\tilde{\boldsymbol{\Phi}}$ as
	\begin{equation}
	\boldsymbol{\Phi}_{i,j,{\ell}}^s = \left\{ \begin{array}{lcc}
	\boldsymbol{\Phi}_{i,j,{\ell}}^s & if \; \frac{1}{NML}\sum_{i,j,{\ell}} \boldsymbol{\Phi}_{i,j,{\ell}}^s < T_{\epsilon} 
	\\ 0 & \mbox{otherwise},
	\end{array}
	\right. 
	\end{equation}
	where $T_{\epsilon}=0.1$ is a customizable hyperparameter that represents the minimum accepted transmittance value for a CA. 

	In this manner, we learn and determine the number of snapshots, as the number of resulting non-all-zero CAs.
	
	Figure \ref{fig:trade} \textcolor{black}{(left)} compares the quality \textcolor{black}{improvement of including the learning of $S$} in terms of PSNR along a spaced interval of resulting compression ratios from $0.1$ to $0.4$, obtained by \textcolor{black}{varying the initial regularization parameter $\rho^0$, against fixing the number of snapshots to be used per each compression ratio,} and training the same network configuration. It can be seen that the most significant improvement of the proposed method occurs at low compression ratios, obtaining a gain of up to $4$ dB while presenting a similar performance for high levels of compression ratios. Further, Fig. \ref{fig:trade} \textcolor{black}{(right)} depicts some of the obtained CAs showing that some of the returned CAs by the E2E approach with the proposed regularizer contain all-zero elements.

		\subsection{{Number of trainable parameters experiment}}
\begin{figure}[t!]
		\begin{tikzpicture}
		 \tikzstyle{every node}=[font=\small]
		 \node[inner sep=0pt] (Static) at (0\linewidth,0){
	 \includegraphics[width=1\linewidth]{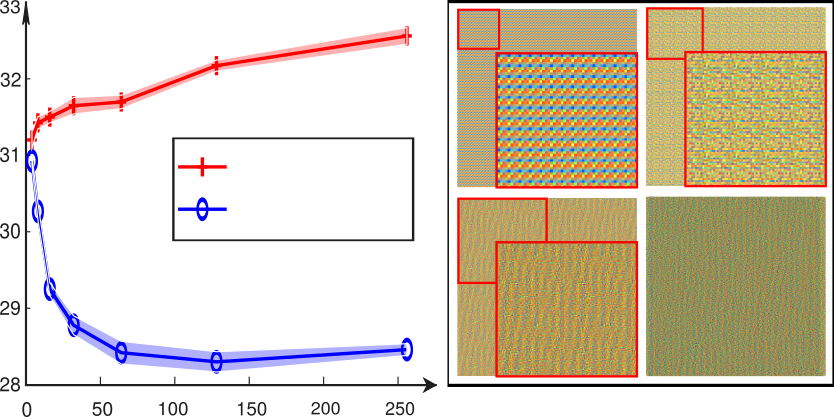}};
	 \node at (-0.21\linewidth, 2.4) {\textbf{Kernel size vs Quality}};
	 \node at (0.25\linewidth, 2.4) {\textbf{Resulted CAs}};
	 
	 \node at (0.025\linewidth,-2.1) { $\mathbf{Q}$};
	 \node at (-0.44\linewidth,2.25) {[dB]};
	 
	 \node[anchor=west] at (-0.22\linewidth,0.4 ) {Training Set};
	 \node[anchor=west] at (-0.22\linewidth, -0.05) {Testing Set};
	 
	 
	 \end{tikzpicture}\vspace{-1.2em}
	 \caption{Effectiveness measured in PSNR for the Training and Testing set, for different number of parameters represented in the kernel size. }
	 \label{fig:kernel_size}\vspace{-1em}
	\end{figure}
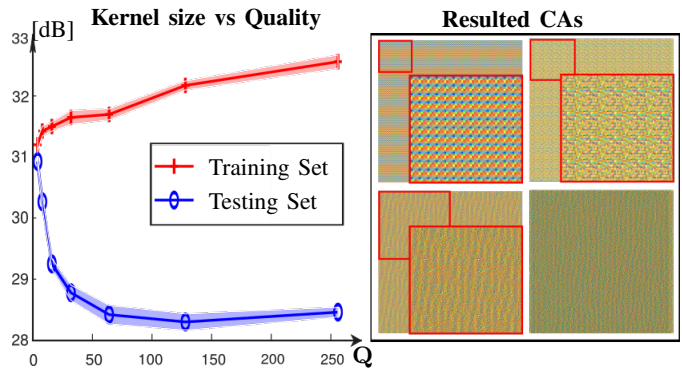

This numerical test aims to demonstrate the importance of reducing the number of trainable parameters to avoid overffiting problems in the E2E approach by using the CA structure strategy presented in \eqref{eq:Kernel_trainable}. For this, we select the CASSI system, with $5<(L=31)$ optical filters in the Colored-CA, for the reconstruction task. Figure. \ref{fig:kernel_size} plots the obtained quality in terms of PSNR when varying the CA kernel size $\mathbf{Q}$ in \eqref{eq:Colored_CA} that directly determines the number of trainable parameters. It can be observed that, as the number of parameters decreases, i.e., the kernel size decreases, the distance between the training and test obtained qualities decreases, indicating a reduced fitting of the model to the training set.


	\subsection{Manufacturing noise experiment}
	
	 	 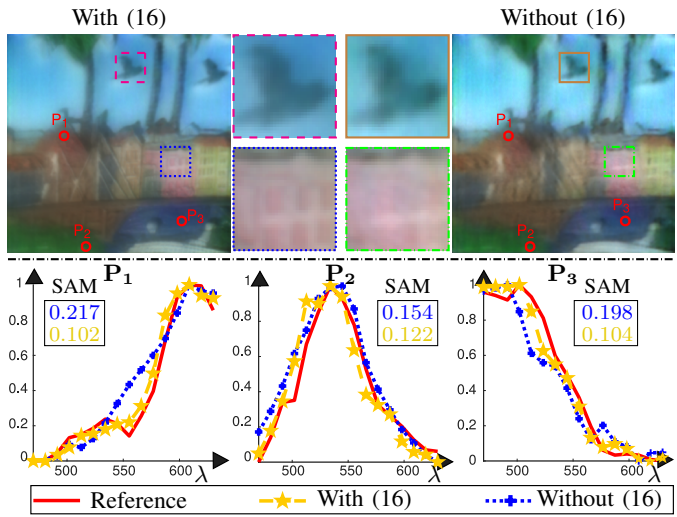
\begin{figure}[t!]
		\begin{tikzpicture}
		 \tikzstyle{every node}=[font=\small]

\begin{scope}[node distance = 1mm, inner sep = 0pt,spy using outlines={rectangle, pink, magnification=3.5, every spy on node/.append style={thick}}]
 \node [align = center](img){\includegraphics[width=1\linewidth]{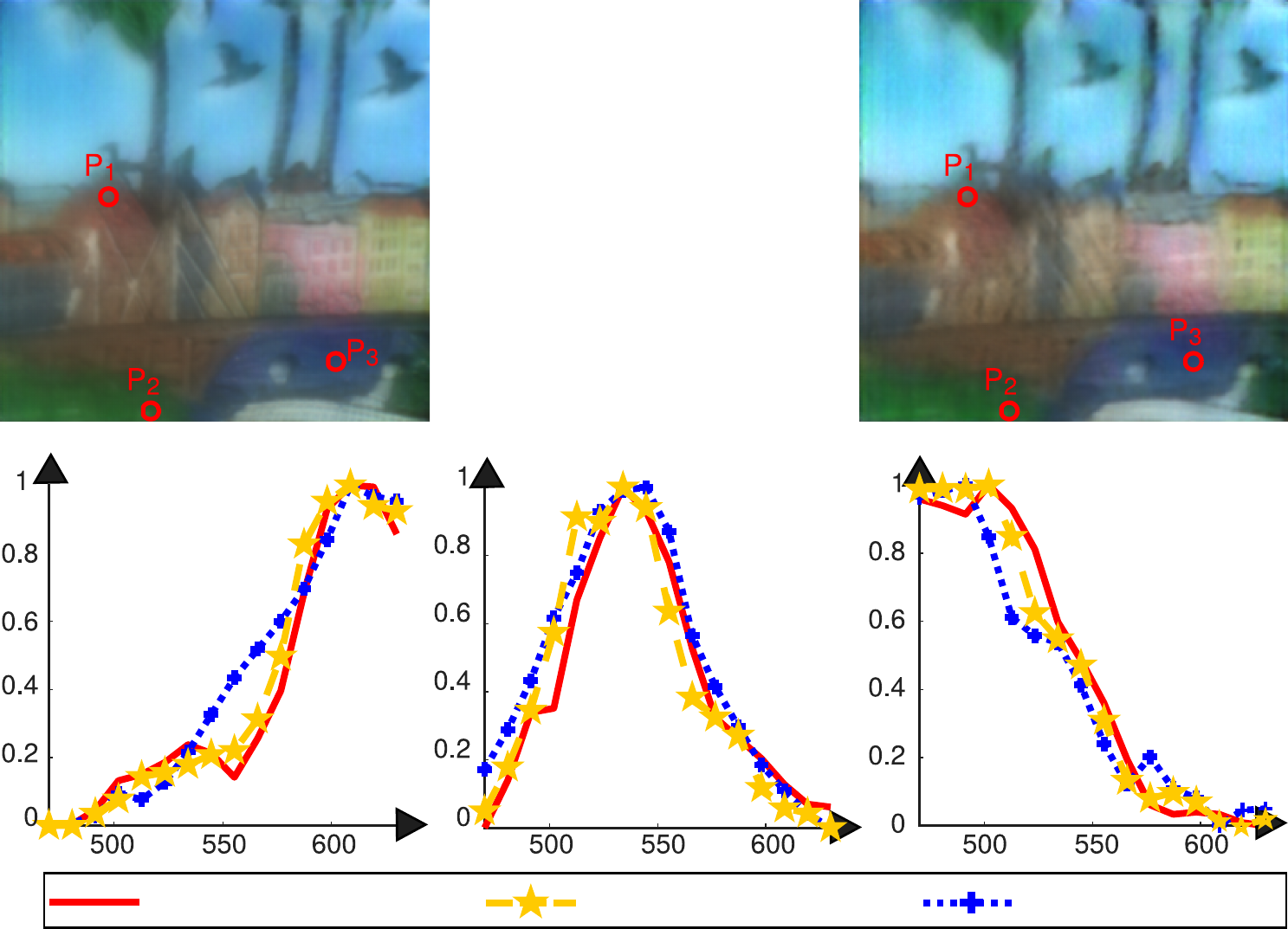}} (img.north west);
 \spy [dashed, magenta, height=0.075\textwidth, width = 0.075\textwidth] on(-2.8,2.75) in node at (-0.75,2.5);
 \spy [solid, brown, thick, height=0.075\textwidth, width = 0.075\textwidth] on(3.1,2.75) in node at (0.75,2.5);
 \spy [blue, densely dotted, thick, height=0.075\textwidth, width = 0.075\textwidth] on(-2.2,1.5) in node at (-0.75,1.0);
 \spy [green, densely dashdotted, thick, height=0.075\textwidth, width = 0.075\textwidth] on(3.7,1.5) in node at (0.75,1.0);
 \end{scope} 
 
 	\node at (-0.3333\linewidth,3.4cm) {With \eqref{eq:ManNois}};
 	\node at (0.3333\linewidth,3.4cm) {Without \eqref{eq:ManNois}};
 	\draw[densely dashdotted,thick] (-0.5\linewidth,0.2) -- (0.5\linewidth,0.2);
		\node at (-0.3333\linewidth,-0.0cm) {$\mathbf{P_1}$};
		\node at (0\linewidth,-0.0cm) {$\mathbf{P_2}$};
	 \node at (0.3333\linewidth,-0.0cm) {$\mathbf{P_3}$};
	 
	 \node at (-0.21\linewidth,-2.65cm) {$\lambda$};
	 \node at (0.14\linewidth,-2.65cm) {$\lambda$};
	 \node at (0.47\linewidth,-2.65cm) {$\lambda$};
	 \node[anchor=west] at (-0.39\linewidth,-3.0cm) {Reference};
	 \node[anchor=west] at (-0.05\linewidth,-3.0cm) {With \eqref{eq:ManNois}};
	 \node[anchor=west] at (0.28\linewidth,-3.0cm) {Without \eqref{eq:ManNois}};

	 \node at (0.1\linewidth,-0.15cm){\footnotesize{SAM}};
\draw (0.055\linewidth,-0.35cm) rectangle (0.145\linewidth,-0.95cm);
\node[blue] at (0.1\linewidth,-0.5cm) {\footnotesize{$0.154$}};
\node[my_yellow] at (0.1\linewidth,-0.8cm) {\footnotesize{$0.122$}};

	 \node at (-0.4\linewidth,-0.15cm){\footnotesize{SAM}};
\draw (-0.445\linewidth,-0.35cm) rectangle (-0.355\linewidth,-0.95cm);
\node[blue] at (-0.4\linewidth,-0.5cm) {\footnotesize{$0.217$}};
\node[my_yellow] at (-0.4\linewidth,-0.8cm) {\footnotesize{$0.102$}};

	 \node at (0.4\linewidth,-0.15cm){\footnotesize{SAM}};
\draw (0.35\linewidth,-0.35cm) rectangle (0.445\linewidth,-0.95cm);
\node[blue] at (0.4\linewidth,-0.5cm) {\footnotesize{$0.198$}};
\node[my_yellow] at (0.4\linewidth,-0.8cm) {\footnotesize{$0.104$}};

		\end{tikzpicture}\vspace{-1.3em}
				\caption{(Top) RGB visual representations of the reconstruction with and without considering the proposed manufacturing noise. (Bottom) Comparison of the normalized spectral signatures at three points in the recovered scenes. At each point it is shown the quantitative SAM metric to quantify the improvement.}
		\label{fig:manufacutre_noise}
		\end{figure}

	 	This experiment evaluates the effectiveness of the manufacturing noise regularizer presented in \eqref{eq:ManNois}. For this, we introduce additive white Gaussian noise with $20$dB of SNR level to perturb the projected encoded measurements, and we introduce uniform random noise between $0$ and $0.2$ to perturb the CA in each forward step for 4 snapshots. Figure \ref{fig:manufacutre_noise} shows a visual comparison of the quality when considering the manufacturing noise in the real testbed setup in Fig. \ref{fig:realsetup}, against not considering the manufacturing noise in the E2E approach. It can be seen a qualitative improvement in the spatial domain, specially in detailed regions of the scene as the shown zoomed sections. Further, in bottom it can be seen the recovered spectral responses at three-points of the scene in which the obtained result when considering the manufacturing noise is closer to the reference one.


	\subsection{Quality improvement via regularizers experiment}
	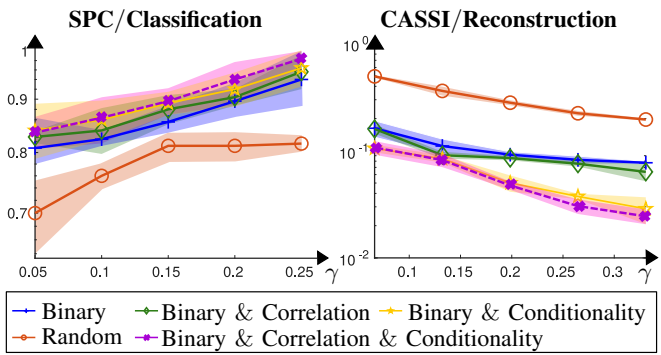
\begin{figure}[b!]
	\begin{tikzpicture}
	\tikzstyle{every node}=[font=\small]
	\node at (0,0.0) {
		\includegraphics[width=0.98\linewidth]{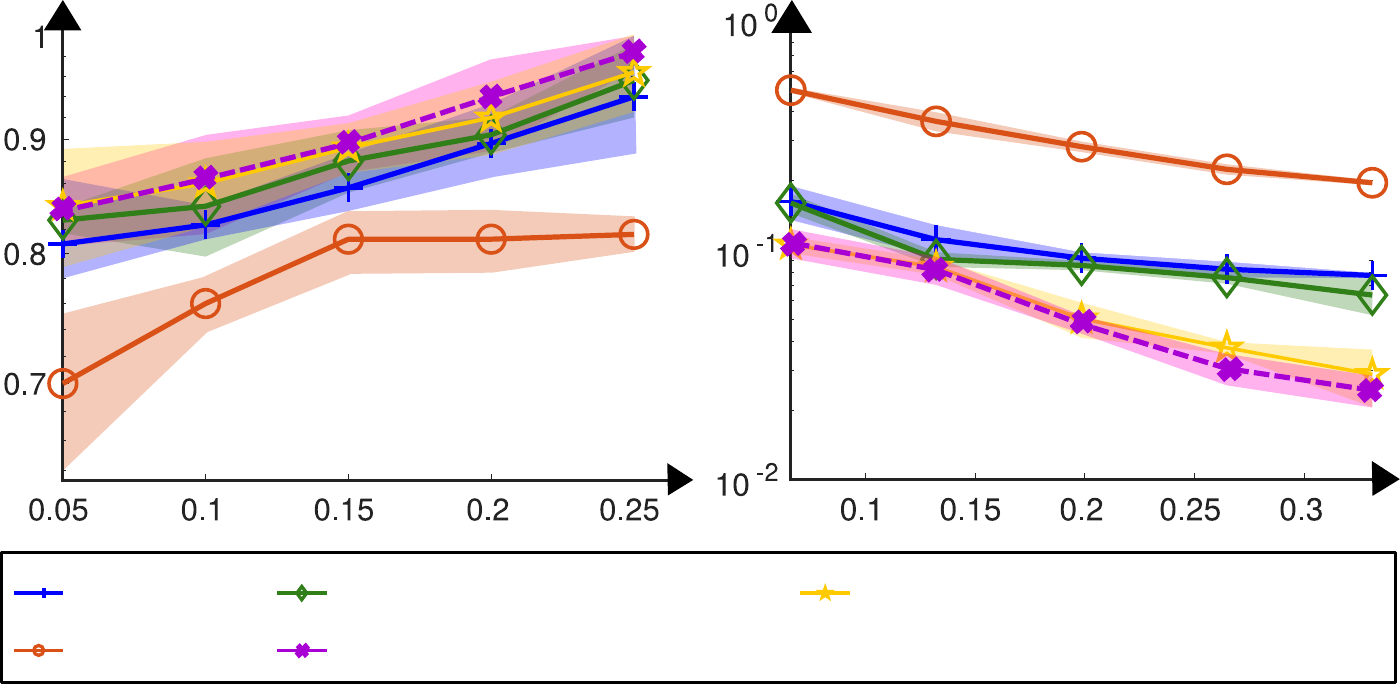}};
	\node at (-0.25\linewidth,2.3cm) {\textbf{SPC$/$Classification}};
	\node at (0.25\linewidth,2.3cm) {\textbf{CASSI$/$Reconstruction}};
	
	\node at (-0.0\linewidth,-1.1cm) {$\gamma$};
	\node at (0.47\linewidth,-1.1cm) {$\gamma$};
	
 	\node[anchor=west] at (-0.45\linewidth,-1.6cm) {Binary};
 	
	\node[anchor=west] at (-0.45\linewidth,-1.9cm) {Random};
	\node[anchor=west] at (-0.27\linewidth,-1.6cm) {Binary $\&$ Correlation};
	\node[anchor=west] at (0.1\linewidth,-1.6cm) {Binary $\&$ Conditionality};
	\node[anchor=west] at (-0.27\linewidth,-1.95cm) {Binary $\&$ Correlation $\&$ Conditionality};

	\end{tikzpicture}\vspace{-1.6em}
		\caption{Quality behavior of adding various regularizers. \textit{(Left)} Real SPC setup for classification task with the quality measured in terms of accuracy. \textit{(Right)} Real CASSI setup for reconstruction task with quality measured in terms of the SAM metric of the spectral response at three spatial locations.}
\label{fig:temporales_3}
\end{figure}	
This experiment evaluates the effectiveness of binary \eqref{eq:family_binary_0}, correlation \eqref{eq:correlation}, and conditionality \eqref{eq:conditionality} regularizers to guide the solution of the optimization problem to a better quality in the CI task for real setups. \textcolor{black}{The initial value of each regularization was set as $\rho _0=1e^{-7}$ except for the binary that was as $\rho _0=1e^{-9}$.} \textcolor{black}{For this, we used the CASSI and a SPC testbed implementation, where $50$ random numbers were printed and used as the target (See supplementary material for further implementation details).} Figure \ref{fig:temporales_3} compares the task quality across an interval of different compression ratios for five cases: using a random CA; using only binary regularizer; using binary and correlation regularizers; using binary and conditionality regularizers; and using binary, correlation, and conditionality regularizers. It can be seen that the quality of the task increases as more regularizers are taken into account, such that using the three regularizers together suppressively outperforms the recovery task quality in both setups. We also remark that the binary regularizer itself provides a significant improvement in comparison to the random scenario.

	\section{Conclusions and Discussions}
	
A customizable E2E approach for CI tasks based on optical coding systems focusing on the CA design has been proposed. The approach enables the consideration of various and extremely important physical and analytical considerations of the CA that can be customized by including specific functions that regularizes the main loss function in the E2E method. The versatility and effectiveness of the proposal was validated over three optical coding systems and three CI tasks covering the customization of features in the CA as binary nature of the elements, ideal transmittance level, ideal number of snapshots, reduced correlation, and improved conditionality. 

\textcolor{black}{We remark the versatility of the proposed regularizations for different applications aiming different purposes, such as the design of binary or real-valued codifications, or capturing different fields of the light.} Furthermore, the proposed E2E scheme can be extended to other optical coding architectures and other sensing systems that include a CA element such as seismic data~\cite{mosher2012compressive} and radar signals~\cite{de2019compressed}; it can also be extended to further applications such as, \textcolor{black}{estimation of depth maps, object detection,} or video processing in which the temporal correlations can be exploited, or to consider the relative position between the CA and the sensor; as well as it can be extended to consider previously formulated regularizers as the uniform sensing presented in~\cite{mejia2018binary, bacca2019super} \textcolor{black}{and the KL-divergence in~\cite{levin2007image}.} Finally, the proposed method is limited to work on real-valued domain, in such a manner that further analysis have to be done to use it on applications such as phase-retrieval in which the complex nature of the CA ensembles conduct to critical differences.

	\ifCLASSOPTIONcaptionsoff
	\newpage
	\fi

	\bibliographystyle{IEEEtran}
	\bibliography{IEEEabrv,Template}
	\vspace{-3em}
	
		\begin{IEEEbiography}[{\vspace{-0.15in}\includegraphics[width=1in,height=1.45in,clip,keepaspectratio]{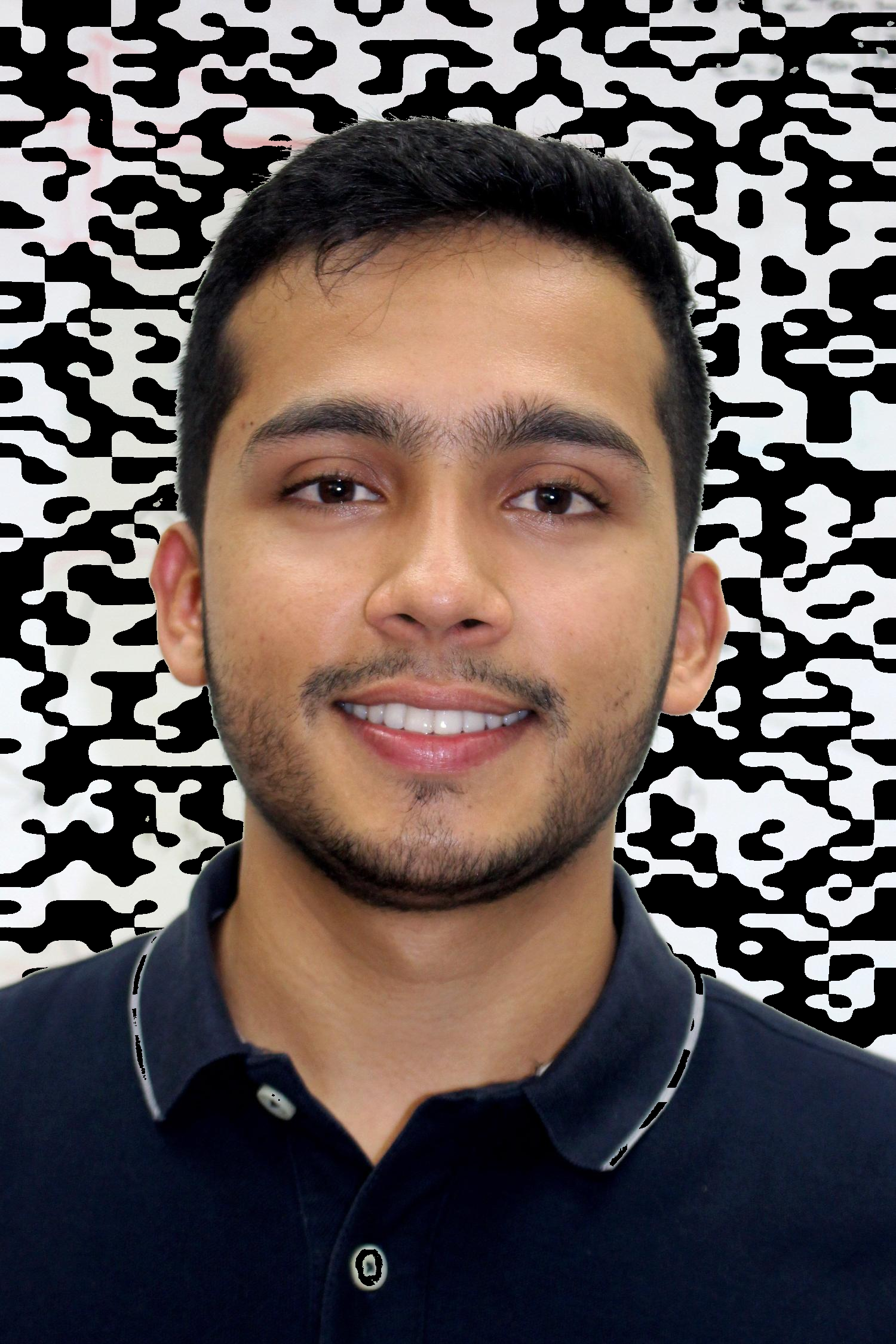}}]{Jorge Bacca}(S'17) received the B.S. degree in
computer science from the Universidad Industrial de Santander, Bucaramanga, Colombia, in 2017, where he is currently pursuing the Ph.D. degree with the Department of Computer Science. His
current research interests include inverse problem, deep learning methods,optical imagining, compressive sensing, phase retrieval, hyperspectral imaging and spectral clustering.
		\end{IEEEbiography}
		\vspace{-3em}
		
	\begin{IEEEbiography}[{\vspace{-0.0in}\includegraphics[width=1in,height=1.45in,clip,keepaspectratio]{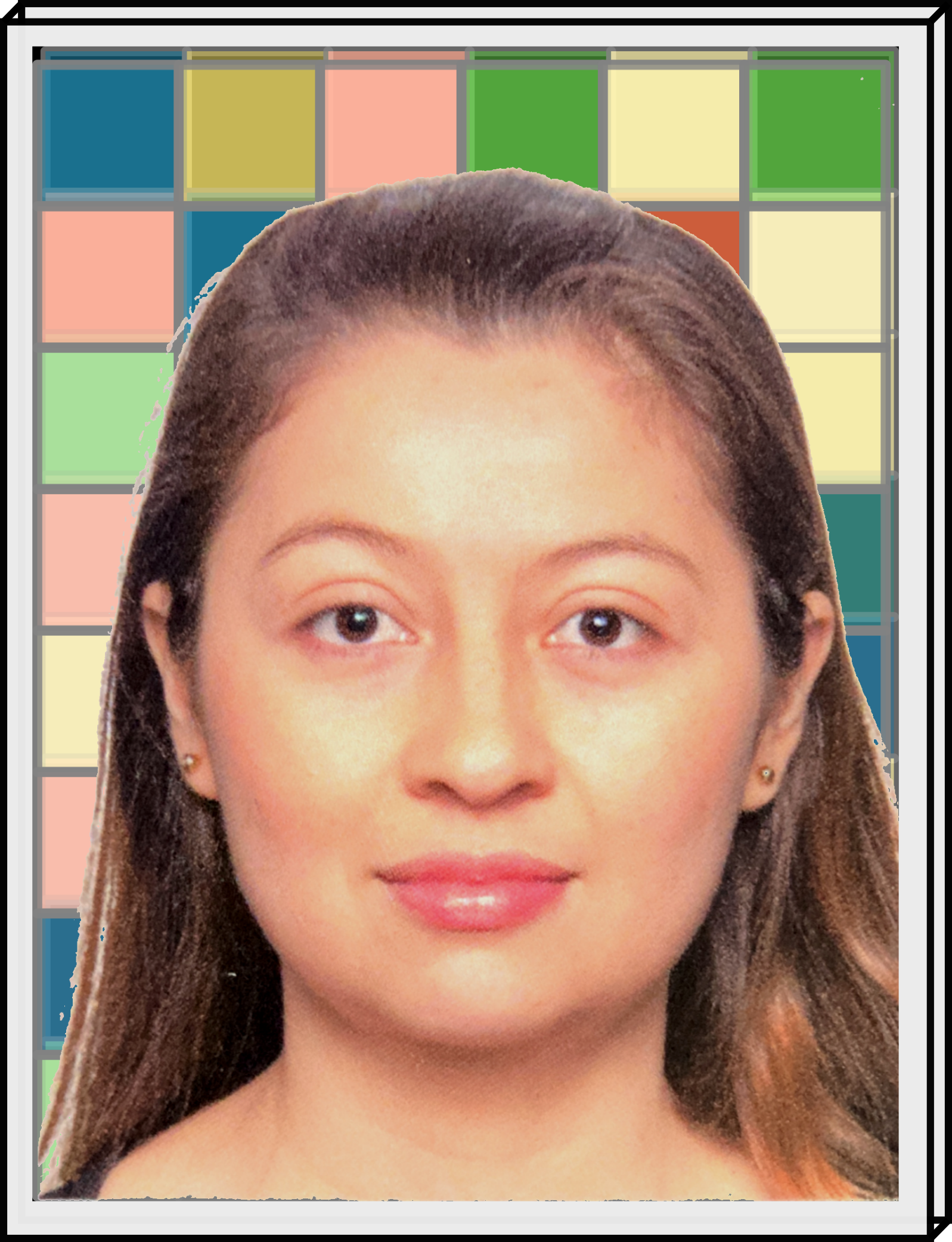}}]{Tatiana Gelvez}
	(S'17) received the B.S. degree in industrial engineering and systems engineering from the Universidad Industrial de Santander, Bucaramanga, Colombia, in 2016, where she is currently pursuing the Ph.D. degree with the Department of Electrical Engineering. During the second semester of 2019 and 2020 she was an intern at Tampere University, Tampere, Finland. Her research interests include numerical optimization, high-dimensional signal processing, spectral imaging, and compressive sensing.
		\end{IEEEbiography}
		\vspace{-3em}
	\begin{IEEEbiography}[{\vspace{-0.05in}\includegraphics[width=1in,height=1.45in,clip,keepaspectratio]{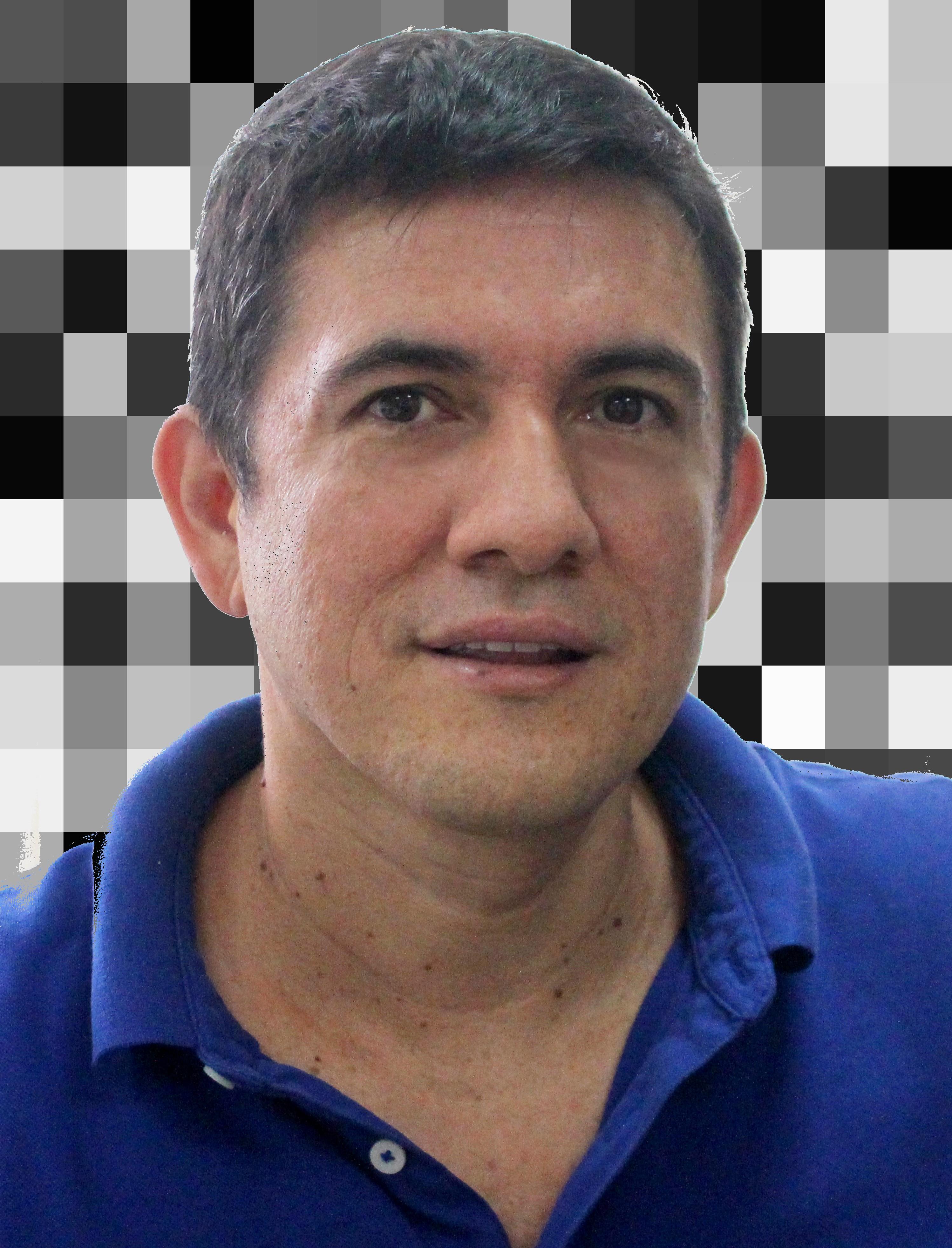}}]{Henry Arguello}
	(S'11–M'13–SM'17) received the B.Sc. Eng. degree in electrical engineering and the M.Sc. degree in electrical power from the Universidad Industrial de Santander, Bucaramanga,
Colombia, in 2000 and 2003, respectively, and the Ph.D. degree in electrical engineering from the University of Delaware, Newark, DE, USA, in 2013. He is currently an Associate Professor with the Department of Systems Engineering, Universidad Industrial de Santander. In first semester 2020, he was a Visiting Professor with Stanford University, Stanford, CA, USA, funded by Fulbright. His research interests include high-dimensional signal processing, optical imaging, compressed sensing, hyperspectral imaging, and CI.
		\end{IEEEbiography}

\end{document}